\documentclass{amsart}

\usepackage{amsmath,amssymb}
\usepackage[latin1]{inputenc}
\usepackage[dvips]{graphics}
\input{xy}
\xyoption{poly}
\xyoption{2cell}
\xyoption{all}
\usepackage{epsfig}  
\usepackage{color}

\newcommand{\ot}{\leftarrow}
\newcommand{\ssi}{\Leftrightarrow}

\newcommand{\za}{\alpha}
\newcommand{\zb}{\beta}
\newcommand{\zd}{\delta}
\newcommand{\zD}{\Delta}
\newcommand{\ze}{\epsilon}
\newcommand{\zg}{\gamma}
\newcommand{\zG}{\Gamma}

\newcommand{\zs}{\sigma}
\newcommand{\zS}{\Sigma}

\newcommand{\Hom}{\textup{Hom}}

\newcommand{\pibar}{\overline{\pi}}

\newtheorem{thm}{Theorem}[section]
\newtheorem{prop}[thm]{Proposition}
\newtheorem{conj}[thm]{Conjecture}
\newtheorem{cor}[thm]{Corollary}
\newtheorem{lem}[thm]{Lemma}
\newtheorem{example}[thm]{Example}
\newtheorem{definition}{Definition} 
\newtheorem{rem}[thm]{Remark}

\newenvironment{pf}{{Proof}.}

\begin{document}
\title{On cluster algebras arising from unpunctured surfaces II}
  \thanks{The author is
   supported by the NSF grant DMS-0700358 and by the University
  of Connecticut} 
\author{Ralf Schiffler}

\date{\today}


\begin{abstract} We study cluster algebras  with principal and arbitrary
  coefficient systems that are associated to unpunctured surfaces. 
We give a direct formula for the Laurent polynomial
expansion of cluster variables in these cluster algebras in terms of
certain paths on a triangulation of the surface.
 As an immediate consequence, we prove the positivity
conjecture of Fomin and Zelevinsky for these cluster algebras.

Furthermore, we obtain direct formulas for $F$-polynomials and $g$-vectors and show that $F$-polynomials have constant term equal to $1$. As an application, we compute the Euler-Poincar\'e characteristic of quiver Grassmannians in Dynkin type $A$ and affine Dynkin type $\tilde A$.
 \end{abstract}
\maketitle



\begin{section}{Introduction}\label{sect intro}

Cluster algebras have been introduced by Fomin and Zelevinsky in \cite{FZ1} in order to create an algebraic framework for total positivity and canonical bases in semisimple algebraic groups. 
Today the theory of cluster algebras is connected to many different areas of mathematics, for example, representation theory of finite dimensional algebras, Lie theory, Poisson geometry and Teichm\"uller Theory.

Cluster algebras are commutative algebras with a distinguished set  of generators, the \emph{cluster variables}. The set of all cluster variables is constructed recursively from an initial set of $n$ cluster variables using so-called \emph{mutations}. Every mutation defines a new cluster variable as a rational function of the cluster variables constructed previously; thus recursively, every cluster variable is a certain rational function in the initial $n$ cluster variables. Fomin and Zelevinsky have shown in \cite{FZ1} that these rational functions are actually Laurent polynomials.

The  first main result of this paper is a direct formula for these Laurent polynomials for the class of cluster algebras that are associated to oriented unpunctured Riemann surfaces with boundary. Let us point out that this formula holds for arbitrary coefficient systems.

In order to be more precise, we need some notation. The cluster variables are grouped into sets of constant cardinality $n$, the \emph{clusters}. The integer $n$ is called the \emph{rank} of the cluster algebra. The cluster algebra is determined by its initial \emph{seed} which consists of a cluster $\mathbf{x}=\{x_1,x_2,\ldots,x_n\}$ together with a \emph{coefficient tuple} $\mathbf{y}=\{y_1,y_2,\ldots,y_n\}$, and a skew-symmetrizable $n\times n$ integer matrix $B=(b_{ij})$. The coefficients $y_1,y_2,\ldots,y_n$ are taken in a torsion free abelian group $\mathbb{P}$. The mutation in direction $k$ defines a new cluster $\mathbf{x}'=\mathbf{x}\setminus\{x_k\}\cup\{x_k'\}$, where

\begin{equation}\label{eq intro}
x_k'=\frac{1}{x_k} \left( y^+\prod_{b_{ki}>0} x_i^{b_{ki}} + y^-\prod_{b_{ki}<0} x_i^{-b_{ki}} \right),
\end{equation}  
where $y^+,y^-$ are certain monomials in the $y_1,y_2,\ldots,y_n$; the precise definition is given in section \ref{sect definition}. Mutations also transform the coefficient tuple $\mathbf{y}$ and the matrix $B$.

If $u$ is any cluster variable, thus $u$ is obtained from the initial cluster $\{x_1,\ldots,x_n\}$ by a sequence of mutations,
then, by \cite{FZ1}, $u$ can be written as a Laurent polynomial in the variables $x_1,x_2,\ldots,x_n$, that is, 
\begin{equation}\label{eq intro2}
 u=\frac{f(x_1,x_2,\ldots,x_n)}{\prod_{i=1}^n x_i^{d_i}},
\end{equation}  
where $f$ is a polynomial with coefficients in the group ring $\mathbb{ZP}$ of the coefficient group $\mathbb{P}$.

Inspired by the work of Fock and Goncharov \cite{FG1,FG2,FG3} and
Gekhtman, Shapiro and Vainshtein \cite{GSV1,GSV2} which discovered
cluster structures in the context of
Teichm\"uller theory, Fomin, Shapiro and Thurston  \cite{FST} and \cite{FT} initiated a
systematic study of the cluster algebras arising from oriented Riemann surfaces with boundary and marked points, a vast class of cluster algebras.

 In this approach, 
cluster variables correspond to isotopyclasses of certain curves in the surface that are called \emph{arcs}.
Clusters are in bijection with triangulations, which are maximal sets of non-crossing arcs, see section \ref{sect surfaces} for precise definitions.
The cluster algebras arising from surfaces all share the property that the number of different matrices that can be obtained from the initial matrix by sequences of mutations is finite \cite[Corollary 12.2]{FST}. Cluster algebras with this property are called \emph{mutation finite}.
While surfaces provide a large class of examples, not all mutation finite cluster algebras are given by surfaces, for example the Dynkin types $E_6,E_7,E_8$ and the affine Dynkin types $\tilde E_6, \tilde E_7, \tilde E_8$ are not.

Assume from now on that the surface has no punctures, which means that all marked points lie on the boundary of the surface.
The cluster algebras arising from these unpunctured surfaces form a three parameter family, the parameters being the genus of the surface, the number of boundary components, and the number of marked points. Fix an arbitrary triangulation $T$ and denote the corresponding cluster by $\mathbf{x}_T$. Let $x_\zg$ be an arbitrary cluster variable, where $\zg$ is an arc in the surface. We know by equation (\ref{eq intro2}) that $x_\zg$ is a Laurent polynomial in $\mathbf{x}_T$.

Following \cite{S,ST,MS}, we introduce \emph{complete $(T,\zg)$-paths} as paths given by concatenation of arcs of the triangulation such that the arcs in even positions are precisely the arcs of $T$ that are crossed by $\zg$ in order, and moreover, between any two such crossings, the complete $(T,\zg)$-path is homotopic to $\zg$.

Suppose for a moment that the cluster algebra has principal coefficients. Then the Laurent polynomial of equation (\ref{eq intro2}) can be written simply as

\[x_\zg =\sum_\za x(\za)\,y(\za),\]
where the sum is over all complete $(T,\zg)$-paths $\za$, where $x(\za)$ is the quotient of the cluster variables associated to the odd arcs in $\za$ by the cluster variables associated to the even arcs in $\za$, and where $y(\za)$ is a product of elements of the initial coefficient tuple that is depending on the orientation of the path $\za$, see Theorem \ref{thm main}.

Using results from \cite{FZ4}, we deduce a formula for the Laurent polynomials in cluster algebras with arbitrary coefficients, see section \ref{sect arbitrary}.
As an immediate corollary, we prove, for cluster algebras  arising from unpunctured surfaces, that Fomin and Zelevinsky's positivity conjecture of 2002 holds, which states that the coefficients in the Laurent polynomials are non-negative integer linear combinations of elements of the coefficient group. This result was shown for finite type cluster algebras in \cite{FZ4}; the only finite type arising from unpunctured surfaces is the Dynkin type $A$, where the corresponding surface is a polygon. In
\cite{CK} and \cite{CR} the conjecture was shown in the case where
the cluster algebra has trivial coefficients and the initial seed is acyclic.

Our formula for the Laurent polynomials in Theorem \ref{thm main} is a continuation of \cite{ST}, where a similar formula was shown for cluster algebras with a very limited coefficient system that was associated to the boundary of the surface. The very special case where the surface is a polygon and coefficients arise from the boundary was covered  in \cite{S}
and also in unpublished work \cite{CP,FZ5} . Our two papers \cite{S,ST} used (ordinary) $(T,\zg)$-paths; the \emph{complete} $(T,\zg)$-paths that we use here where introduced in \cite{MS}. The paper \cite{MS} gives a new parametrization of the formulas in terms of perfect matchings.

The proof of Theorem \ref{thm main} uses the theory of covering spaces and involves cluster algebras with a specially designed coefficient system which combines  principal coefficients and boundary coefficients.

A second main result of this article consists in explicit formulas for the $F$-polynomials and the $g$-vectors defined in \cite{FZ4}, see Theorem \ref{thm F-polynomial} and Theorem \ref{thm g-vectors}. The formula for the $F$-polynomial yields a proof for a conjecture of \cite{FZ4} stating that $F$-polynomials have constant term $1$.

Theorem \ref{thm main} has interesting intersections with work of other
people. In \cite{CCS2}, the authors obtained a formula for the
denominators of the cluster expansion in types $A,D$ and $E$, see also
\cite{BMR}. In \cite{CC,CK,CK2} an expansion formula was given in the
case where the cluster algebra is acyclic and the cluster lies in an
acyclic seed. Palu generalized this formula to arbitrary
clusters in an acyclic cluster algebra \cite{Palu}.  These formulas
use the cluster category introduced in \cite{BMRRT}, and in \cite{CCS1} for
type $A$, and do not give information about the coefficients. 

Recently, Fu and Keller generalized this formula further to cluster algebras  with principal coefficients that admit a categorification by a 2-Calabi-Yau category \cite{FK}, and, combining results of \cite{A} and \cite{ABCP,LF}, such a categorification exists in the case of cluster algebras associated to unpunctured surfaces. We thank Bernhard Keller for pointing out this last fact to us. The formula of Keller and Fu involve the Euler-Poincar\'e characteristic of quiver Grassmannians.
Comparing it to our result, we derive a formula for this Euler-Poincar\'e characteristic  in Dynkin type $A$ and affine Dynkin type $\tilde A$, see Theorem \ref{thm EPC}. In particular, we show that the Euler-Poincar\'e characteristic is non-negative, see Corollary \ref{cor EPC}. We also study projective presentations in the module categories of finite dimensional algebras associated to the triangulations of the surfaces.

In \cite{SZ,CZ,Z,MP} cluster expansions for cluster algebras of
rank 2 are given. 
A related cluster expansion formula in the type 
$A$ case was given in \cite{Propp}.  

In \cite{M} a cluster expansion for cluster algebras of
finite classical type is given for clusters that lie in a bipartite seed.

The paper is organized as follows. In section \ref{sect cluster algebras}, we recall some definitions and results from the theory of cluster algebras. We state our formula for the Laurent polynomials in cluster algebras with principal coefficients in Theorem \ref{thm main} in section \ref{sect expansion}. In that section, we also introduce complete $(T,\zg)$-paths.
Sections \ref{sect proof} and \ref{sect cover} are devoted to the proof of Theorem \ref{thm main}; section \ref{sect proof} deals with the simply connected case and section \ref{sect cover} with the general case. Section \ref{sect cover} also contains some necessary preparatory results on Galois coverings of quivers.
In section \ref{sect Fg}, we present our formulas for $F$-polynomials and $g$-vectors, and prove Conjectures 5.4 and 5.5 of \cite {FZ4}.
A formula for the Laurent expansions in cluster algebras with arbitrary coefficients is given in section \ref{sect arbitrary}. In that section, we also prove the positivity conjecture for arbitrary coefficients of geometric type.
In section \ref{sect EPC}, we study Euler-Poincar\'e characteristics of quiver Grassmannians and projective presentations in the module categories of finite dimensional algebras associated to triangulated surfaces in \cite{ABCP}.

The author thanks Gregg Musiker and Hugh Thomas for several interesting discussions.

\end{section}

\begin{section}{Cluster algebras}\label{sect cluster algebras}
In this section, we recall some facts on cluster algebras. For further
details, the reader is referred to \cite{FZ4,FST}.

\begin{subsection}{Definition}\label{sect definition}
Following Fomin and Zelevinsky, let $(\mathbb{P},\oplus,\cdot)$ be a \emph{semifield}, i.e.
 an abelian group $(\mathbb{P},\cdot)$ together with a binary
 operation $\oplus$ which is commutative, associative, and, 
 with respect to the multiplication $\cdot$, 
 distributive.
 Let $\mathbb{ZP}$ be the group ring of $(\mathbb{P},\cdot)$ and let
 $\mathcal{F}=\mathbb{QP}(x_1,x_2,\ldots,x_n)$ be the field of
 rational functions in $n$ variables with coefficients in
 $\mathbb{QP}$.

 A \emph{seed} $\zS$ is a triple $\zS=(\mathbf{x},\mathbf{y},B)$,
 where
\begin{itemize}
\item[-] $\mathbf{x}=\{x_1,x_2,\ldots,x_n\}$ is a transcendence basis
  of $\mathcal{F}$ over $\mathbb{QP}$,
\item[-]  $\mathbf{y}=\{y_1,y_2,\ldots,y_n\}$ is an $n$-tuple of
  elements $y_i\in \mathbb{P}$, and
\item[-] $B=(b_{ij})$ is a skew symmetrizable $n\times n$ integer matrix.
\end{itemize}   

The set $\mathbf{x}$ is called a \emph{cluster} and its elements are
\emph{cluster variables}. The set $\mathbf{y}$ is called
\emph{coefficient tuple} and $B$ is called \emph{exchange matrix}.

\begin{rem} In this paper we will study a special type of cluster
  algebras, those arising from surfaces, and for these cluster
  algebras  the matrix $B$ is always
  skew symmetric.
\end{rem}

Given a seed $ (\mathbf{x},\mathbf{y},B)$ its mutation
$\mu_k(\mathbf{x},\mathbf{y},B)$ in direction $k$ is a new seed
$(\mathbf{x}',\mathbf{y}',B')$ defined as follows. Let
$[x]_+=\textup{max}(x,0)$. 

\begin{itemize}
\item[-] $B'=(b_{ij}')$ with 
\begin{equation}\label{eq 90}
 b_{ij}'=\left\{
\begin{array}{ll}
-b_{ij} &\textup{if $i=k$ or $j=k$,}\\
b_{ij}+[-b_{ik}]_+\,b_{kj}+b_{ik}[b_{kj}]_+ &\textup{otherwise,}
\end{array}  \right.
\end{equation}
\item[-] $\mathbf{y}'=(y_1',y_2',\ldots,y_n')$ with 
\begin{equation}\label{eq 91} y_j'=\left\{\begin{array}{ll}
y_k^{-1} &\textup{if $j=k$,}\\
y_jy_k^{[b_{kj}]_+}(y_k\oplus 1)^{-b_{kj}}&\textup{if $j\ne k$,}
\end{array}\right.
\end{equation}
\item[-] $\mathbf{x}'=\mathbf{x}\setminus \{x_k\}\cup\{x_k'\}$ where
\begin{equation}\label{exchange relation}
  x_k'= \frac{y_k\prod x_i^{[b_{ik}]_+} + \prod
  x_i^{[-b_{ik}]_+}}{(y_k\oplus 1)\, x_k}.
\end{equation} 
\end{itemize}        
The formula (\ref{exchange relation}) is called \emph{exchange
  relation}. 
One can check that the mutations are involutions, that is, $\mu_k\mu_k
(\mathbf{x},\mathbf{y},B)=(\mathbf{x},\mathbf{y},B)$. 

Most of the time, we will be dealing with cluster algebras of
\emph{geometric type}, which means that $\mathbb{P}$ is a tropical
semifield $\mathbb{P}=\textup{Trop}(u_1,\ldots,u_\ell)$, 
that is, $(\mathbb{P},\cdot)$ is a free abelian group on the
  generators 
$u_{1},u_{2},\ldots,u_{\ell}$, and  the addition $\oplus$ is given by
  the formula
\[\prod_{j=1}^\ell u_j^{a_j}\oplus\prod_{j=1}^\ell u_j^{b_j} 
= \prod_{j=1}^\ell u_j^{\textup{min}(a_j,b_j)}.
\]
In cluster algebras of geometric type, it is convenient to replace the
matrix $B$ by an $(n+\ell)\times n$ matrix $\tilde B=(b_{ij})$ whose upper part
is the $n\times n$ matrix $B$ and whose lower part is an $\ell\times
n$ matrix that encodes the coefficient tuple via

\begin{equation}\label{eq 20}
y_k = \prod_{i=1}^{\ell} u_i^{b_{(n+i)k}}.
\end{equation}  
Then the mutation of the coefficient tuple in equation (\ref{eq 91}) is determined by the mutation
of the matrix $\tilde B$ in equation (\ref{eq 90}) and the formula (\ref{eq 20}); and the
exchange relation (\ref{exchange relation}) becomes
\begin{equation}\label{geometric exchange}
 x_k'=x_k^{-1} \left( \prod_{i=1}^n x_i^{[b_{ik}]_+}
\prod_{i=1}^{\ell} u_i^{[b_{(b+i)k}]_+} 
+\prod_{i=1}^n x_i^{[-b_{ik}]_+}
\prod_{i=1}^{\ell} u_i^{[-b_{(n+i)k}]_+}
\right).
\end{equation}  
If the cluster algebra is of geometric type then the group ring
 $\mathbb{ZP}$ of $(\mathbb{P},\cdot)$,  
 is the ring of Laurent polynomials in the variables
$u_1, u_2, \ldots, u_\ell$.

Two seeds $\zS_1,\zS_2$ are called \emph{mutation equivalent} if there is a sequence
of mutations $\mu=\mu_{i_1}\mu_{i_2}\cdots\mu_{i_s}$ such that
$\mu\zS_1=\zS_2$. Thus starting from an initial seed
$(\mathbf{x},\mathbf{y},B)$ one constructs the class of all seeds that
are mutation equivalent to the initial one by successive mutations in
all possible directions.

Define $\mathcal{X}$ to be the set of all cluster variables, that is,
$\mathcal{X}$ is the union of all clusters $\mathbf{x}'$
such that 
there exists a seed $(\mathbf{x}',\mathbf{y}',B')$ that is mutation
equivalent to the initial seed.  The cluster algebra
$\mathcal{A}(\mathbf{x},\mathbf{y},B)$ is the $\mathbb{ZP}$-subalgebra
of the field $\mathcal{F}$ generated by the set of all cluster
variables, thus 
\[\mathcal{A}(\mathbf{x},\mathbf{y},B)=\mathbb{ZP}[\mathcal{X}] .\]

Using the exchange relations, each cluster variable can be written as
a rational function of the cluster variables in the initial seed. The
following theorem, known as the Laurent phenomenon, states that these
rational functions are actually Laurent polynomials.

\begin{thm}\label{Laurent phenomenon}\cite[Theorem 3.1]{FZ1}
Let $x\in \mathcal{X} $ be any cluster variable in the cluster algebra
$\mathcal{A}(\mathbf{x},\mathbf{y},B)$. Then $x$ has an expansion in 
the initial cluster $\mathbf{x}=\{x_1,x_2,\ldots,x_n\}$ as
\[ x=\frac{f(x_1,x_2,\ldots,x_n)} {x_1^{d_1} x_2^{d_2}\cdots
  x_n^{d_n}}, 
\]
where the right hand side is a reduced fraction and $f\in
\mathbb{ZP}[x_1,x_2,\ldots,x_n], d_i\ge 0$.
\end{thm}  

It has been conjectured in \cite{FZ1} that the polynomials $f$ have
non-negative coefficients.
\begin{conj} [Positivity Conjecture] Each coefficient of the
  polynomial $f$ in Theorem \ref{Laurent phenomenon}
  is a non-negative integer linear combination of elements in $\mathbb{P}$.
\end{conj}

\end{subsection} 

\begin{subsection}{Cluster algebras with principal
    coefficients}\label{sect principal coefficients}
Fomin and Zelevinsky introduced in \cite{FZ4} a special type of
coefficients, called \emph{principal coefficients}. We recall some of
their properties here.

A cluster algebra $\mathcal{A}=\mathcal{A}(\mathbf{x},\mathbf{y},B)$ is
said to have \emph{principal coefficients} if
the coefficient semifield $\mathbb{P}$ is the tropical semifield
$(\mathbb{P},\oplus,\cdot)=\textup{Trop}(y_{1},y_{2},\ldots,y_{n})$
with the initial coefficient tuple
$\mathbf{y}=\{y_{1},y_{2},\ldots,y_{n}\}$ as set of generators.

In particular, $\mathcal{A}$ is of geometric type, so one can replace
the $n\times n$ matrix $B$ by a $2n\times n$ matrix $\tilde B$ whose
upper part is the $n\times n $ matrix $B$ and whose lower part encodes 
the coefficient tuple $\mathbf{y}$ in the seed by formula (\ref{eq
  20}). In this description, the cluster algebra $\mathcal{A}$ has
principal coefficients if the lower part of the initial matrix $\tilde
B$ is the $n\times n$ identity matrix.

\begin{prop}\label{prop principal laurent}
 For cluster algebras with principal coefficients, the cluster
expansions of Theorem \ref{Laurent phenomenon}  take the form
\[x=\frac{f(x_1,x_2,\ldots,x_n;y_1,y_2,\ldots,y_n)}{x_1^{d_1},x_2^{d_2},\ldots,x_n^{d_n}},
\] 
where $f\in\mathbb{Z}[x_1,x_2\ldots,x_n;y_1,y_2,\ldots,y_n]$.
\end{prop} 

\begin{pf} \cite[Proposition 3.6]{FZ4}
\qed
\end{pf}  

Knowing the cluster expansions for a cluster algebra with principal
coefficients allows one to compute the cluster expansions for the
``same'' cluster algebra with an arbitrary coefficient system. More
precisely, let $\hat A=\mathcal{A}(\mathbf{x},\hat{\mathbf{y}},B)$ be
a cluster algebra over $\mathbb{Z}\hat{\mathbb{P}}$ with initial
coefficient tuple $\hat{\mathbf{y}}=(\hat 
y_1,\hat y_2,\ldots,\hat y_n)$, where
$\hat{\mathbb{P}}$ is an arbitrary 
semifield. Denote by $\hat{\mathcal{F}}$ the 
field of rational functions in $n$ variables with coefficients in
$\mathbb{Q}\hat{\mathbb{P}}$. Let
$A=\mathcal{A}(\mathbf{x},\mathbf{y},B)$ be the cluster algebra with
principal coefficients that has the  same initial cluster and the same
initial exchange matrix as $\hat {\mathcal{A}}$. Let $\hat x $ be any
cluster 
variable in $\hat{\mathcal{A}}$ and let
$\mu=\mu_{i_1}\mu_{i_2}\cdots\mu_{i_s}$ be a sequence of mutations
such that $\hat x\in \mu(\mathbf{x})$ but $\hat
x\notin\mu_{i_2}\cdots\mu_{i_s}(\mathbf{x})$. 
Using the same sequence of mutations in $\mathcal{A}$, let $x$ be the unique cluster variable in $\mathcal{A}$
such that $ x\in \mu(\mathbf{x})$ but $
x\notin\mu_{i_2}\cdots\mu_{i_s}(\mathbf{x})$, and let
\[x=\frac{f(x_1,x_2,\ldots,x_n;y_1,y_2,\ldots,y_n)}{x_1^{d_1}\ldots x_n^{d_n}},
\] 
 be the cluster   
expansion of $x$ in the initial seed $(\mathbf{x},\mathbf{y},B) $.

\begin{thm}\label{thm FZ4}\cite[Theorem 3.7]{FZ4} With the above notation, the
  cluster expansion of $\hat x$ in the initial cluster
  $(\mathbf{x},\hat{\mathbf{y}},B) $ in $\hat{\mathcal{A}}$ is 
\[\hat x=\frac
{f\vert_{\hat{\mathcal{F}}}(x_1,x_2,\ldots,x_n;\hat {y}_1,\hat
  {y}_2,\ldots,\hat {y}_n)}
{x_1^{d_1}\ldots x_n^{d_n}\,
  f\vert_{\hat{\mathbb{P}}}(1,1,\ldots,1;\hat y_1,\hat
  y_2,\ldots,\hat y_n)},
\]
where $f\vert_{\hat{\mathcal{F}}}(x_1,x_2,\ldots,x_n;\hat y_1,\hat
  y_2,\ldots,\hat y_n)$ is the polynomial $f$ evaluated in the field
  $(\hat{\mathcal{F}},+,\cdot)$ after substituting $\hat y_i$ for
  $y_i$, 
  and $ f\vert_{\hat{\mathbb{P}}}(1,1,\ldots,1;\hat y_1,\hat
  y_2, \ldots,\hat y_n)$ is the polynomial $f$ evaluated in the
  semifield 
  $(\hat{\mathbb{P}},\oplus,\cdot)$ after substituting $1$ for $x_i$
  and 
  $\hat y_i$ for $y_i$, $i=1,2,\ldots,n$.
\end{thm}  
Thus the cluster expansion of $\hat x$ is obtained from the expansion
for $x$ by replacing $y_i$ by $\hat y_i$, and dividing by  $
f\vert_{\hat{\mathbb{P}}}(1,1,\ldots,1;\hat y_1,\hat 
  y_2, \ldots,\hat y_n)$.
\end{subsection} 

\begin{subsection}{Cluster algebras arising from unpunctured
    surfaces}\label{sect surfaces} 

In this section, we recall the construction of \cite{FST} in the
case of surfaces without punctures.

Let $S$ be a connected oriented 2-dimensional Riemann surface with
boundary and $M$ a non-empty set of marked points in the closure of
$S$ with at least one marked point on each boundary component. The
pair $(S,M)$ is called \emph{bordered surface with marked points}. Marked
points in the interior of $S$ are called \emph{punctures}.  

In this paper we will only consider surfaces $(S,M)$ such that all
marked points lie on the boundary of $S$, and we will refer to $(S,M)$
simply by \emph{unpunctured surface}. The orientation of the surface will play a crucial role. 

We say that two curves in $S$ \emph{do not cross} if they do not intersect
each other except that endpoints may coincide.

\begin{definition}
An \emph{arc} $\zg$ in $(S,M)$ is a curve in $S$ such that 
\begin{itemize}
\item[(a)] the endpoints are in $M$,
\item[(b)] $\zg$ does not cross itself,
\item[(c)] the relative interior of $\zg$ is disjoint from $M$ and
  from the boundary of $S$,
\item[(d)] $\zg$ does not cut out a monogon or a digon. 
\end{itemize}   
\end{definition}     
 Curves that connect two
marked points and lie entirely on the boundary of $S$ without passing
through a third marked point are called \emph{boundary arcs}.
Hence an arc is a curve between two marked points, which does not
intersect itself nor the boundary except possibly at its endpoints and
which is not homotopic to a point or a boundary arc.

Each arc is considered up to isotopy inside the class of such curves. Moreover, each arc is considered up to orientation, so if an arc has endpoints $a,b\in M$ then it can be represented by a curve that runs from $a$ to $b$, as well as by a curve that runs from $b$ to $a$.

For any two arcs $\zg,\zg'$ in $S$, let $e(\zg,\zg')$ be the minimal
number of crossings of $\zg$ and $\zg'$, that is, $e(\zg,\zg')$ is the
minimum of
the numbers of crossings of  arcs $\za$ and $\za'$, where $\za$ is
isotopic to $\zg$ and $\za'$ is isotopic to $\zg'$.
Two arcs $\zg,\zg'$ are called \emph{compatible} if $e(\zg,\zg')=0$. 
A \emph{triangulation} is a maximal collection of
compatible arcs together with all boundary arcs. 
The arcs of a 
triangulation cut the surface into \emph{triangles}.
Since $(S,M)$ is an unpunctured surface, the three sides of each
triangle are distinct (in contrast to the case of surfaces with
punctures).  Any triangulation  has
$n+m$ elements, $n$ of which  are arcs in $S$, and the remaining $m$
elements are boundary arcs. Note that the number of boundary arcs
is equal to the number of marked points.

\begin{prop}\label{prop rank}
The number $n$ of arcs in any triangulation is  given by the formula 
$n=6g+3b+m-6$,  where $g$ is the
genus of $S$, $b$ is the number of boundary components and $m=|M|$ is the
number of marked points. The number $n$ is called the \emph{rank} of $(S,M)$.
\end{prop}  
\begin{pf} \cite[2.10]{FST}
\qed
\end{pf}  

\begin{cor}\label{cor triangles}
The number of triangles in any triangulation is equal to 
\[ n-2(g-1)-b.
\]
\end{cor}
\begin{pf}
Consider the Riemann surface without boundary obtained from $S$ by gluing a disk into each boundary component. 
Computing the Euler-Poincar\'e characteristic of this surface using the genus on the one hand and using the simplicial complex given by the triangulation $T$ on the other hand leads to the equation
\[ 2-2g = \textup{number of triangles} +b -(n+m) +m,
\] 
and the statement follows.
\qed
\end{pf}
Note that $b> 0$ since the set $M$ is not empty.
 Table \ref{table 1} gives some examples of unpunctured surfaces.

\begin{table}
\begin{center}
  \begin{tabular}{ c | c | c || l  }
  \  b\ \  &\ \  g \ \   & \ \  m \ \  &\  surface \\ \hline
    1 & 0 & n+3 & \ polygon \\ 
    1 & 1 & n-3 & \ torus with disk removed \\
    1 & 2 & n-9 & \ genus 2 surface with disk removed \\\hline 
    2 & 0 & n & \ annulus\\
    2 & 1 & n-6 & \ torus with 2 disks removed \\ 
    2 & 2 & n-12 & \ genus 2 surface with 2 disks removed \\ \hline
    3 & 0 & n-3 & \ pair of pants \\ \\
  \end{tabular}
\end{center}
\caption{Examples of unpunctured surfaces}\label{table 1}
\end{table}

Following  \cite{FST}, we associate a cluster algebra
to the unpunctured surface $(S,M)$ as follows.
 Choose any triangulation
$T$, let $\tau_1,\tau_2,\ldots,\tau_n$ be the $n$ interior arcs of
$T$ and  denote the  $m$ boundary
arcs of the surface by $\tau_{n+1},\tau_{n+2},\ldots,\tau_{n+m}$. 
For any triangle $\Delta$ in $T$ define a matrix 
$B^\Delta=(b^\Delta_{ij})_{1\le i\le n, 1\le j\le n}$  by
\[ b_{ij}^\Delta=\left\{
\begin{array}{ll}
1 & \textup{if $\tau_i$ and $\tau_j$ are sides of 
  $\Delta$ with  $\tau_j$ following $\tau_i$  in the }
\\ &\textup{ counter-clockwise order;}\\
-1 &  \textup{if $\tau_i$ and $\tau_j$ are sides of
  $\Delta$ with  $\tau_j$ following $\tau_i$  in the }\\
&\textup{  clockwise  order;}\\
0& \textup{otherwise.}
\end{array} \right. \]
Then define the matrix 
$ B_{T}=(b_{ij})_{1\le i\le n, 1\le j\le n}$  by
$b_{ij}=\sum_\Delta b_{ij}^\Delta$, where the sum is taken over all
triangles in $T$. Note that the boundary arcs of the triangulation are ignored in the definition of $B_{T}$.
Let $\tilde B_{T}=(b_{ij})_{1\le i\le 2n, 1\le j\le n}$ be the
$2n\times n$ matrix whose upper $n\times n$ part is $B_{T}$ and whose
lower $n\times n$ part is the identity matrix.
The matrix $B_{T}$ is skew-symmetric and each of its entries  $b_{ij}$ is either
$0,1,-1,2$, or $-2$, since every arc $\tau$ can be in at most two triangles. 
An example where $b_{ij}=2 $ is given in Figure
\ref{fig bij=2}.

\begin{figure}
\centering
\input{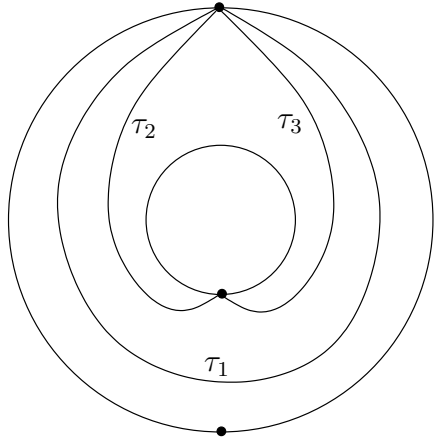}
\caption{A triangulation with $b_{23}=2$ \label{fig bij=2}}
\end{figure}

Let $\mathcal{A}(\mathbf{x}_{T},\mathbf{y}_{T},B_{T})$ be the cluster algebra with principal
coefficients in the triangulation $T$, that is,  $\mathcal{A}(\mathbf{x}_{T},\mathbf{y}_{T},B_{T})$ is
given by the seed
$(\mathbf{x}_{T},\mathbf{y}_{T},B_{T})$ where
$\mathbf{x}_{T}=\{x_{\tau_1},x_{\tau_2},\ldots,x_{\tau_n}\}$ is the
cluster associated to the triangulation $T$, and the initial coefficient
vector $\mathbf{y}_{T}=(y_1,y_2,\ldots,y_n)$ is the vector of generators
of $\mathbb{P}=\textup{Trop}(y_1,y_2,\ldots,y_n)$. 

For the  boundary arcs we define $x_{\tau_k}=1$,  $k=n+1,n+2,\ldots,n+m$.

For each $k=1,2,\ldots,n$,  there is a unique quadrilateral  in $T\setminus \{\tau_k\}$
in which $\tau_k$ is one of the diagonals. Let $\tau_k'$ denote the
other diagonal in that quadrilateral.
Define the
\emph{flip} $\mu_k T$ to be the triangulation
$T\setminus\{\tau_k\}\cup\{\tau_k'\}$.
The mutation $\mu_k$ of the seed $\zS_T$ in the cluster algebra
$\mathcal{A}$ corresponds to the flip $\mu_k$ of the triangulation $T$
in the following sense.
The matrix $\mu_k(B_T)$ is the matrix corresponding to the
triangulation $\mu_k T$,
the cluster
$\mu_k(\mathbf{x}_T)$ is $\mathbf{x}_T\setminus\{x_{\tau_k}\}\cup
\{x_{\tau_k'}\}$, and the corresponding exchange relation is given by
\[x_{\tau_k} x_{\tau_k'} = x_{\rho_1} x_{\rho_2} y^+ +  x_{\zs_1}
x_{\zs_2} y^-,
\]
where $y^+,y^-$ are some coefficients, and $\rho_1,\zs_1,\rho_2,\zs_2$
are the sides of the quadrilateral in which $\tau_k$ and $\tau_k'$ are
the diagonals, such that  $\rho_1,\rho_2$ are opposite sides and
$\zs_1,\zs_2$ are opposite sides too.

\end{subsection}

\begin{subsection}{The quiver of a triangulation}\label{sect quiver}

Since the cluster algebra $\mathcal{A}=\mathcal{A}(\mathbf{x}_{T},\mathbf{y}_{T},B_{T})$ defined in section \ref{sect surfaces} has principal coefficients, we can consider the $(2n\times n)$ matrix $\tilde B_{T}$ whose upper part is $B_{T}$ and whose lower part is the $n\times n$ identity matrix. If $\mu$ is a sequence of mutations, then the upper part of $\mu \tilde B_{T}$ corresponds to the triangulation $\mu T$, but the triangulation 
 does not give any information about the
lower $n\times n$ part of the matrix $\tilde B_T$.

In order to keep track of the whole matrix and, thus, the coefficient tuple,
we find it convenient to use quivers. 

Let $T$ be the triangulation that will serve as initial seed for the
cluster algebra with principal coefficients. Define a quiver $Q_T$ as
follows. The vertices of $Q_T$ are labeled by integers
$1,2,\ldots,2n$, where the first $n$ vertices correspond to the
interior arcs $\tau_1,\ldots,\tau_n$ of $T$, and the second $n$
vertices correspond to the initial coefficient tuple $y_1,\ldots,y_n$.
The arrows of $Q_T$ are given by the matrix $\tilde B_T=(\tilde
b_{ij})$, that is, for each pair $i>j$ of vertices, the quiver $Q_T$
has $\tilde b_{ij}$ arrows from $i$ to $j$; 
where we use the convention that, if $\tilde b_{ij}<0$ then we have
$-\tilde b_{ij} $ arrows from $j$ to $i$.

In \cite{ABCP}, generalizing a
construction of \cite{CCS1}, the 
authors define relations for the full subquiver of 
$Q_T$ whose set of vertices is 
$\{1,2,\ldots,n\}$ and show that the corresponding bound quiver
algebra is a gentle algebra. In the case where the surface is a
polygon, respectively an annulus, these algebras are precisely the
cluster-tilted algebras of type $A_n$, respectively $\tilde A_n$. In
section \ref{sect projective presentations}, we will study projective
presentations of certain indecomposable modules over these algebras.

Independently, these relations have also been defined in \cite{LF}
using quivers with potentials introduced in \cite{DWZ}. In \cite{LF},
the author considers the more general situation where the surface is
allowed to have punctures. In the punctured case, the resulting bound
quiver algebras are no longer gentle.
\end{subsection} 

\end{section} 

\begin{section}{Cluster expansions with principal
  coefficients}\label{sect expansion} 
In this section we state our formula for the cluster expansions with principal coefficients in Theorem \ref{thm main}.

Let $T=\{\tau_1,\ldots,\tau_n,\tau_{n+1},\ldots,\tau_{n+m}\}$ be a
triangulation of the unpunctured surface $(S,M)$, where
$\tau_1,\ldots,\tau_n$ are arcs and $\tau_{n+1},\ldots,\tau_{n+m}$ are
boundary arcs ($m=|M|$). Let $B_T$ be the corresponding $2n\times n$
matrix with lower half equal to the $n\times n$ identity matrix, and
let $\mathcal{A}=\mathcal{A}(\mathbf{x}_T,\mathbf{y}_T,B_T)$ be the
cluster algebra  with principal coefficients and initial seed
$(\mathbf{x}_T,\mathbf{y}_T,B_T)$, where
$\mathbf{x}_T=\{x_{\tau_1},\ldots,x_{\tau_n}\}$ is the
  initial cluster and 
$\mathbf{y}_T=\{y_{\tau_1},\ldots,y_{\tau_n}\} $ is the initial  coefficient vector.
We will often write $x_i$ and $y_i$ instead of $x_{\tau_i}$ and
$y_{\tau_i}$ respectively. By \cite{FST}, the cluster variables in $\mathcal{A}$
correspond to the arcs in $(S,M)$. 

Let $\zg$ be any arc in $(S,M)$ that crosses $T$ exactly $d$ times. We
fix an orientation for $\zg$ and we denote its starting point by $s$
and its endpoint by $t$, with $s,t\in M$, see Figure \ref{fig notation}. Let
$s=p_0,p_1,\ldots,p_d,p_{d+1}=t$ 
be the intersection  points of $\zg$ and $T$ in order of occurrence on
$\zg$, hence $p_0, p_{d+1}\in M$ and each $p_i$ with $1\le i\le d$
lies in the interior of $S$.
 Let $i_1,i_2,\ldots,i_d$ be such that $p_k$ lies on the arc
$\tau_{i_k}\in T$, for $k=1,2,\ldots,d$. Note that $i_k$ may be equal
to $i_j$ even if $k\ne j$.

For $k=0,1,\ldots,d$, let $\zg_k$ denote the segment of the path $\zg$
from  the point $p_k$ to the point $p_{k+1}$. Each $\zg_k$ lies in
exactly one triangle $\zD_k$ in $T$. If $1\le k\le d-1$, the triangle
$\zD_k$ is formed by the arcs $\tau_{i_k},
\tau_{i_{k+1}}$ and a third arc that we denote by
$\tau_{[\zg_k]}$. In the triangle $\zD_{0}$,  $\tau_{i_1}$ is one of the sides. Denote
the side of $\zD_{0}$ that lies clockwise of  $\tau_{i_1}$ by  ${\tau_{[\zg_0]}}$ and the
other side by  $\tau_{[\zg_{-1}]}$.
Similarly,  $\tau_{i_d}$ is one of the sides of $\zD_d$. Denote
the side that lies clockwise of  $\tau_{i_d}$ by  $\tau_{[\zg_{d}]}$ and the
other side by  $\tau_{[\zg_{d+1}]}$.
\begin{figure}[htbp]
\begin{center}
\input{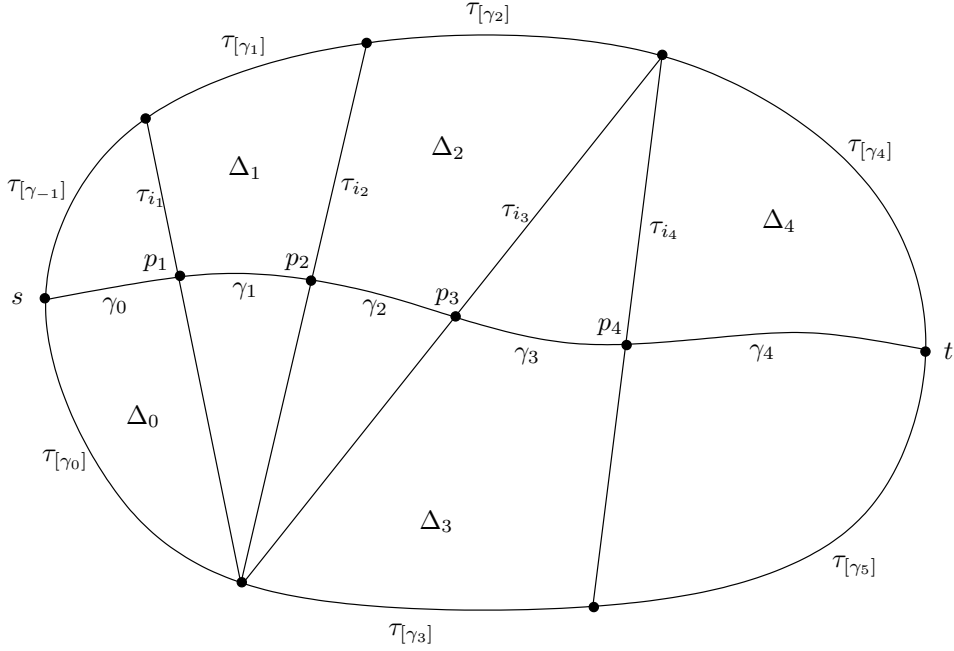}
\caption{Relative position of $\zg$ and $T$, an example with $d=4$}
\label{fig notation}
\end{center}
\end{figure}

\begin{subsection}{Complete $(T,\zg)$-paths}\label{sect Tg}

A \emph{$T$-path} is a path $\za$ in $S$ on the triangulation $T$,
that is, there exist arcs $\za_1,\za_2,\ldots,\za_{\ell(\za)}\in T$ such
that $\za$ is the concatenation of paths
$\za=\za_1\za_2\cdots\za_{\ell(\za)}$. We will write
$\za=(\za_1,\za_2,\ldots,\za_{\ell(\za)})$.
Recall that arcs in $(S,M)$ are not oriented. Now a $T$-path induces an orientation on each of its arcs $\za_i$. Note that a $T$-path $\za$ may not be uniquely determined by the sequence $(\za_1,\za_2,\ldots,\za_{\ell(\za)})$ since arcs may be loops and, thus, picking a starting point and an endpoint does not determine the orientation. However, the sequences $(\za_1,\za_2,\ldots,\za_{\ell(\za)})$ determine the \emph{complete} $(T,\zg)$-paths that we are going to define now.
\begin{definition}\cite{MS}
A $T$-path $\alpha=(\alpha_1,\alpha_2,\dots ,
\alpha_{\ell(\alpha)}) $ is called  a \emph{complete $(T,\gamma)$-path}
if the following axioms hold: 

\begin{itemize}
\item[(T1)] The even arcs are precisely the arcs crossed by $\zg$ in
  order, that is, $\za_{2k}=\tau_{i_k}$.
\item[(T2)] For all $k=0,1,2,\ldots,d$, the segment $\zg_k$ is
  homotopic to the segment of the path $\za$ starting at the point
  $p_k$ following $\za_{2k},\za_{2k+1}$ and $\za_{2k+2}$ until the
  point $p_{k+1}$.
\end{itemize}

\end{definition}
 \begin{rem}\begin{itemize}
\item Every complete $(T,\zg)$-path starts and ends at the same point
  as $\zg$, because of (T2).
\item Every complete $(T,\zg)$-path has length $2d+1$.
\item  For all arcs $\tau$ in the triangulation $T$, the number of
  times that
$\tau$ occurs as $\alpha_{2k}$ is exactly the number
  of crossings between $\zg$ and $\tau$.
\item In contrast to the ordinary $(T,\zg)$-paths introduced in \cite{ST}, complete $(T,\zg)$-paths allow backtracking.
\end{itemize}
\end{rem} 

\end{subsection}

\begin{subsection}{Orientation}\label{sect orientation}
By Corollary \ref{cor triangles}, the triangulation $T$ cuts the surface $S$ into $(n-2(g-1)-b)$ triangles. The orientation of the surface $S$ induces an orientation on each of these triangles in such a way that, whenever two triangles $\zD,\zD'$ share an
edge $\tau$, then the orientation of $\tau $ in $\zD$ is opposite to
the orientation of $\tau $ in $\zD'$, see Figure \ref{figdisc}. 
There are precisely two such orientations, we assume without loss of generality that we have the ``clockwise
orientation'', that is, in each triangle $\zD$, going around the
boundary of $\zD$ according to the orientation of $\zD$ is clockwise
when looking at it from outside the surface.

\begin{figure}
\includegraphics{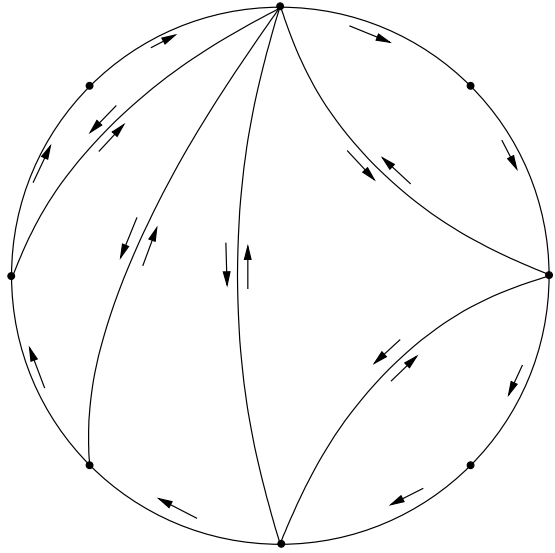}
\caption{A triangulation of the disk with $8$ marked points with
  clockwise orientation}\label{figdisc}
\end{figure}

Let $\za $ be a complete $(T,\zg)$-path. Then
$\za_{2k}=\tau_{i_k}$ is a common edge of the two triangles
$\zD_{k-1}$ and $\zD_k$.
We say that $\za_{2k}$ is $\zg$-\emph{oriented} if the orientation of
$\za_{2k}$ in the path $\za$ is the same as the orientation of
$\tau_{i_k}$  in the triangle $\zD_k$, see Figure \ref{fig
  gammaoriented}.
\begin{figure}
\input{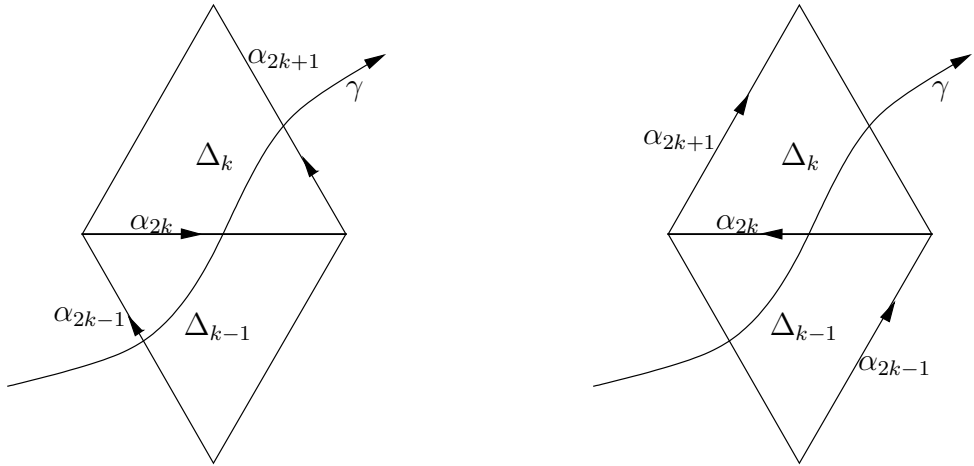}
\caption{Two examples of the  $(T,\zg)$-path segment
  $(\za_{2k-1},\za_{2k},\za_{2k+1})$. On the left, 
  $\za_{2k}$ is not $\zg$-oriented and on the right, $\za_{2k}$ is $\zg$-oriented.}\label{fig gammaoriented}
\end{figure}   

\end{subsection}

\begin{subsection}{Expansion Formula}\label{sect EF}
We are ready to state the main result of this section. 
We keep the setup of the previous sections.
Define 
\[x(\za) =\frac{\prod_{k \textup{ odd}} x_{\za_k} }{x_{i_1}
  x_{i_2}\ldots x_{i_d}} \qquad \textup{and}\qquad y(\za)= \prod_{k:
  \za_{2k}\textup{ is $\zg$-oriented}} y_{i_k}\]
Thus $x(\za)$ does not depend on the orientation, but $y(\za)$
  does. Recall that $x_{\za_k}=1$ if $\za_k$ is a boundary arc.

 The following theorem
gives a cluster expansion formula of an arbitrary cluster variable
$x_\zg$ in the initial cluster with principal coefficients.

\begin{thm}\label{thm main}
Let $x_\zg $ be any cluster variable in the cluster algebra $\mathcal{A}$. Then its expansion in the initial seed $(\mathbf{x}_T,\mathbf{y}_T,B_T)$ is given by
\[ x_{\zg} = \sum_{\za}x(\za)\ y(\za) \]
where the sum is over all complete $(T,\zg)$-paths $\za$ in $(S,M)$.
\end{thm}  

\begin{rem}
The formula in Theorem \ref{thm main} does not depend on our choice of orientation on the arc $\zg$. Indeed, considering the opposite orientation $\zg^{op}$, there is a bijection between the set of  complete $(T,\zg)$-paths and complete $(T,\zg^{op})$-paths, sending a path $\za=(\za_1,\za_2,\ldots,\za_{2d+1})$ to the opposite path $\za^{op}=(\za_{2d+1},\ldots,\za_{2},\za_1)$; moreover $x(\za^{op})=x(\za)$, and $y(\za^{op})=y(\za)$.
\end{rem}

The proof of Theorem \ref{thm main} will be given in section \ref{sect proof} for simply connected surfaces and in section \ref{sect cover} for arbitrary surfaces. To
illustrate the statement, we give two examples here.

\begin{example}\label{ex A}{\bf The case $A_n$:} The cluster algebra
  $\mathcal{A}$ is of type $A_n$ if $(S,M)$ is an
  $(n+3)$-gon. Our example illustrates the case $n=5$.
The following figure shows a triangulation $T=\{\tau_1,\ldots,\tau_{13}\}$ and
a (dotted) arc $\zg$. Next to it is the list of complete $(T,\zg)$-paths. 
\[\begin{array}{cc}
\def\alphanum{\ifcase \xypolynode \or 6 \or 7\or 8 \or 9\or 10\or
  11\or 12\or 13\fi}
\xy/r6pc/: {\xypolygon8"A"{~<<{@{}}~><{@{-}|@{}}
~>>{_{\tau_{\alphanum}}}}},
\POS"A2" \ar@{-}|(0.7)@{}^(0.7){\tau_1} "A4",
\POS"A4" \ar@{-}|@{}_(0.5){\tau_2} "A6",
\POS"A6" \ar@{-}|(0.3)@{}^(0.3){\tau_3} "A2",
\POS"A2" \ar@{-}|@{}^{\tau_4} "A8",
\POS"A6" \ar@{-}|(0.3)@{}^(0.3){\tau_5} "A8",
\POS"A3" \ar@{.}|(0.6)@{}^(0.6){\zg} "A7",
\endxy
&

\begin{array}{c}
 ( \tau_8,\tau_1,\tau_3,\tau_3,\tau_3,\tau_5,\tau_{12})\\
\\( \tau_8,\tau_1,\tau_3,\tau_3,\tau_4,\tau_5,\tau_{11})\\
\\( \tau_7,\tau_1,\tau_2,\tau_3,\tau_3,\tau_5,\tau_{12}) \\
\\( \tau_7,\tau_1,\tau_2,\tau_3,\tau_4,\tau_5,\tau_{11}) \\
\\( \tau_7,\tau_1,\tau_1,\tau_3,\tau_5,\tau_5,\tau_{11}) \\
\end{array} 
\end{array}  \]
Theorem \ref{thm main} thus implies that 
\[\begin{array}{rcll}
 x_\zg&=& 
({x_1x_3x_5})^{-1}&
( x_3^2\,y_1y_3y_5 \\
&&
+x_3x_4\,y_1y_3\\ &&
+x_2x_3\,y_3y_5\\ &&
+x_2x_4\,y_3\\ &&
+x_1x_5)
\end{array}\]
\end{example}

\begin{example}\label{ex Atilde}{\bf The case $\tilde A_{n-1}$:}
 The cluster algebra
  $\mathcal{A}$ is of type $\tilde A_{n-1}$ if $(S,M)$ is an
  annulus. Our example illustrates the case $n=4$. Figure
  \ref{figaffine} shows a triangulation
  $T=\{\tau_1,\tau_2,\ldots,\tau_8\}$ and a (dotted) arc $\zg$.
The complete list of $(T,\zg)$-paths is as follows:
\[
\begin{array}{ll}
(\tau_4,\tau_1,\tau_2,\tau_2,\tau_2,\tau_3,\tau_4,\tau_4,\tau_4,\tau_1,\tau_2)\\
(\tau_4,\tau_1,\tau_2,\tau_2,\tau_2,\tau_3,\tau_4,\tau_4,\tau_5,\tau_1,\tau_8)\\
(\tau_4,\tau_1,\tau_2,\tau_2,\tau_6,\tau_3,\tau_7,\tau_4,\tau_4,\tau_1,\tau_2)\\
(\tau_4,\tau_1,\tau_2,\tau_2,\tau_6,\tau_3,\tau_7,\tau_4,\tau_5,\tau_1,\tau_8)\\
(\tau_4,\tau_1,\tau_2,\tau_2,\tau_6,\tau_3,\tau_3,\tau_4,\tau_1,\tau_1,\tau_8)\\
\\
(\tau_5,\tau_1,\tau_8,\tau_2,\tau_2,\tau_3,\tau_4,\tau_4,\tau_4,\tau_1,\tau_2)\\
 (\tau_5,\tau_1,\tau_8,\tau_2,\tau_2,\tau_3,\tau_4,\tau_4,\tau_5,\tau_1,\tau_8)\\
 (\tau_5,\tau_1,\tau_8,\tau_2,\tau_6,\tau_3,\tau_7,\tau_4,\tau_4,\tau_1,\tau_2)\\
 (\tau_5,\tau_1,\tau_8,\tau_2,\tau_6,\tau_3,\tau_7,\tau_4,\tau_5,\tau_1,\tau_8)\\
 (\tau_5,\tau_1,\tau_8,\tau_2,\tau_6,\tau_3,\tau_3,\tau_4,\tau_1,\tau_1,\tau_8)\\
\\
 (\tau_5,\tau_1,\tau_1,\tau_2,\tau_3,\tau_3,\tau_7,\tau_4,\tau_4,\tau_1,\tau_2)\\
 (\tau_5,\tau_1,\tau_1,\tau_2,\tau_3,\tau_3,\tau_7,\tau_4,\tau_5,\tau_1,\tau_8)\\
 (\tau_5,\tau_1,\tau_1,\tau_2,\tau_3,\tau_3,\tau_3,\tau_4,\tau_1,\tau_1,\tau_8)\\
\end{array}  
\]
Hence Theorem \ref{thm main} implies
\[
\begin{array}{rcll}
x_\zg &=& \displaystyle (x_1x_2x_3x_4x_1)^{-1} 
(x_4x_2x_2x_4x_4x_2\,y_1y_2y_3y_4y_1\\&&
+x_4x_2x_2x_4  \,y_1y_2y_3y_4\\&&
+x_4x_2  x_4x_2\,y_1y_2y_4y_1\\&&
+x_4x_2    \,y_1y_2y_4\\&&
+x_4x_2 x_3x_1 \,y_1y_2\\&&

\\&&
+  x_2x_4x_4x_2\,y_2y_3y_4y_1\\&&
+  x_2x_4  \,y_2y_3y_4\\&&
+    x_4x_2\,y_2y_4y_1\\&&
+      \,y_2y_4\\&&
+   x_3x_1 \,y_2\\&&

\\&&
+ x_1x_3 x_4x_2\,y_4y_1\\&&
+ x_1x_3   \,y_4\\&&
 x_1x_3x_3x_1 )
\end{array}  
\]
Note that the term $x_2^2x_4^2y_1y_2y_3y_4$ appears with multiplicity two.

\begin{figure}
\begin{center}
\input{figaffine.pstex_t}
\caption{The case $\tilde A_{n-1}$}\label{figaffine}
\end{center}
\end{figure}   
\end{example}

\end{subsection}

\begin{subsection}{Positivity}\label{sect positivity}

The following positivity
conjecture of \cite{FZ1} is a direct 
consequence of Theorem \ref{thm main}.

\begin{cor}\label{cor pos} Let $(S,M) $ be an
  unpunctured surface and let $\mathcal A= 
\mathcal{A}(\mathbf{x},\mathbf{y},B)$ be the cluster algebra with principal
 coefficients in the seed $(\mathbf{x},\mathbf{y},B)$ associated to
 some triangulation of $(S,M)$. Let $u$ be any cluster variable 
 and let
\[x=\frac{f(x_1,\ldots,x_n,y_{1},\ldots,y_{n})}{x_1^{d_1}\ldots x_n^{d_n}}\] be the
expansion of $u$ in the cluster $\mathbf{x}\{x_1,\ldots,x_n\}$, where $f$ is a
polynomial which is not divisible by any of the $x_1,\ldots,x_n$.
Then  the coefficients of $f$ are non-negative integers.
\end{cor}  

\begin{pf} This is a direct consequence of Theorem \ref{thm main}.
\qed
\end{pf}  

\begin{rem}
We will also prove this conjecture for arbitrary coefficients of geometric type in Theorem \ref{thm pos}.
\end{rem}
\end{subsection}

\end{section}




\begin{section}{Proof of Theorem \ref{thm main} for simply connected surfaces}\label{sect proof}
In this section, we prove Theorem \ref{thm main} in the case where the surface
$S$ is simply connected. 

Recall that $T=\{\tau_1,\ldots,\tau_n,\tau_{n+1},\ldots,\tau_{n+m}\}$ is a
triangulation, $\zg$ is  an arc 
on which we fixed an orientation such that $\zg$ is going from $s$ to $t$,
$\zg$ crosses $T$ exactly $d$ times,
$s=p_0,p_1,\ldots,p_d,p_{d+1}=t$ 
are the intersection  points of $\zg$ and $T$ in order of occurrence on
$\zg$, and $i_1,i_2,\ldots,i_d$ are such that $p_k$ lies on the arc
$\tau_{i_k}\in T$, for $k=1,2,\ldots,d$. 

Also recall that $\zg_k$ denotes the segment of the path $\zg$
from  the point $p_k$ to the point $p_{k+1}$. Each $\zg_k$ lies in
exactly one triangle $\zD_k$ in $T$. If $1\le k\le d-1$, the triangle
$\zD_k$ is formed by the arcs $\tau_{i_k},
\tau_{i_{k+1}}$ and a third arc that we denote by
$\tau_{[\zg_k]}$.

\begin{lem}\label{lem 0}
If $S$ is simply connected then
\begin{itemize}
\item[(a)] Any two arcs in $S$ cross at most once.
\item[(b)] Two paths are homotopic if they have the same startpoint
  and the same endpoint. In particular, any arc is uniquely determined
  by its endpoints. 
\item[(c)] $i_k=i_j$ if and only if $k=j$.
\item[(d)] $\zD_k=\zD_j$ if and only if $k=j$.
\end{itemize}    \qed
\end{lem}

In the triangle  $\zD_0$, $\tau_{i_1}$ is one of the sides. Denote
the side of $\zD_0$ that lies clockwise of  $\tau_{i_1}$ by  ${\tau_{[\zg_0]}}$ and the
other side by  $\tau_{[\zg_{-1}]}$, see Figure \ref{fig st}.
Similarly,  $\tau_{i_d}$ is one of the sides of $\zD_d$. Denote
the side that lies clockwise of  $\tau_{i_d}$ by  $\tau_{[\zg_{d}]}$ and the
other side by  $\tau_{[\zg_{d+1}]}$.
\begin{figure}
\input{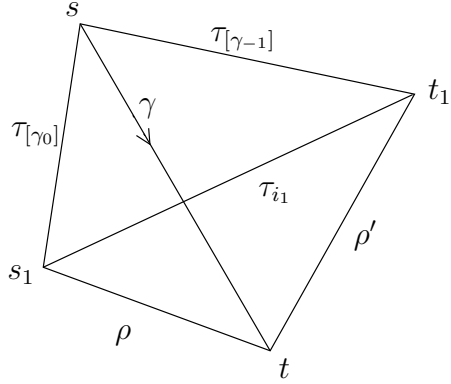}
\caption{Proof of the simply connected case.}\label{fig st}
\end{figure}   
Note that ${\tau_{[\zg_{-1}]}} $ and ${\tau_{[\zg_0]}}$ may be boundary arcs. Let $s_1$  be the
common endpoint of $\tau_{i_1}$ and ${\tau_{[\zg_0]}}$ and let $t_1$ be the
common endpoint of $\tau_{i_1 }$ and ${\tau_{[\zg_{-1}]}}$. Since $S$ is simply
connected, there is a unique arc $\rho$ from $s_1 $ to $t$ and a
unique arc $\rho'$ from $t_1$ to $t$, and, moreover, the arcs
${\tau_{[\zg_{-1}]}},{\tau_{[\zg_0]}},\rho$, and $\rho'$ form a quadrilateral in which $\zg$ and
$\tau_{i_1}$ are the diagonals, see Figure \ref{fig st}.

With this notation, we have the following Lemma.
\begin{lem}\label{lem 5}
If $S$ is simply connected,
then 
any arc $\tau\in T$ that crosses $\rho$ (respectively $\rho'$) also crosses $\zg$.
\end{lem}
\begin{pf}
If $\tau$ crosses $\rho$ then $\tau$ meets the triangle whose sides are $\rho, \tau_{i_1}$ and $\rho'$. But $\tau$ does not cross $\tau_{i_1}$  since both are arcs of the triangulation $T$. Thus $\tau$ must cross $\zg$. The proof of the second statement is similar.
\qed
\end{pf}

Let $\mu_j$ denote the flip in direction $j$.
\begin{lem}\label{lem 1}
If $S$ is simply connected,
then $\zg\in \mu_{i_1}\mu_{i_2}\ldots\mu_{i_d} (T)$.
\end{lem}  

\begin{pf} We prove the statement by induction on $d$. If $d=0$ then
  $\zg\in T$, and there
  is nothing to prove. Suppose $d>0$, then
\[ \mu_{i_d}(T)=T\setminus\{\tau_{i_d}\} \cup \{\tau_{i_d}'\},\]
and $\tau_{i_d}'$ is incident to $t$, hence $\tau_{i_d}'$ does not
cross $\zg$, because $S$ is simply connected. Thus  $\mu_{i_d}(T)$ is  a triangulation and
$p_1,p_2,\ldots,p_{d-1}$ are its crossing points with $\zg$ in
order. By induction, $\zg\in \mu_{i_1}\ldots\mu_{i_{d-1}}\mu_{i_d}(T)$.
\qed
\end{pf}  
 \begin{cor}\label{cor 1}
Let $S$ be simply connected.
Then the arcs $\tau_{i_1},{\tau_{[\zg_{-1}]}},{\tau_{[\zg_0]}},\rho,$ and $\rho'$ are all elements of
$\mu_{i_2}\ldots\mu_{i_d}(T)$. 
\end{cor}  

\begin{pf}
By Lemma \ref{lem 1}, we have $\zg\in \mu_{i_1}\mu_{i_2}\ldots\mu_{i_d}(T)$. On the other hand,
 we have $\zg\notin \mu_{i_2}\mu_{i_3}\ldots\mu_{i_d}(T)$ since
$\tau_{i_1}\in \mu_{i_2}\mu_{i_3}\ldots\mu_{i_d}(T)$ and $\zg$ crosses
$\tau_{i_1}$.  Hence
\[\mu_{i_1}\mu_{i_2}\ldots\mu_{i_d}(T)=\mu_{i_2}\ldots\mu_{i_d}(T)
\setminus\{\tau_{i_1}\}\cup\{\zg\}\]
and thus $\mu_{i_2}\mu_{i_3}\ldots\mu_{i_d}(T)$ contains the arcs
involved in the corresponding exchange relation, namely ${\tau_{[\zg_{-1}]}},{\tau_{[\zg_0]}},\rho,\rho'$.
\qed
\end{pf}  

We shall need some more notation.
If $\rho\notin T$ then let $\ell$ be the least integer such that
$\tau_{i_\ell}$ crosses 
$\rho$, and if $\rho\in T$ then let $\ell=d+1$. Lemma \ref{lem 5} guaranties that the integer $\ell$ is well-defined.   Similarly, if
$\rho'\notin T$ then let $\ell'$ be the least integer such that $\tau_{i_{\ell'}}$
crosses $\rho'$, and if $\rho'\in T$ then let $\ell'=d+1$, see Figure
\ref{fig ell}. Note that if $d=1$ then $\ell=\ell'=2$, and if $d>1$
then one of $\ell, \ell'$ is equal to $2$ and the 
other is greater than $2$. 
Also note that the quiver $Q_T$ defined in section \ref{sect quiver} contains the subquiver  
\[ 
\begin{array}{rcl}
i_1\ot i_2\ot \ldots \ot i_{\ell'-1}       &\textup{if}& \ell=2\\
  i_1\to i_2\to\ldots\to i_{\ell-1}  &\textup{if} & \ell'=2.
\end{array}\]
Moreover, every $\tau_{i_j}$ with $j\ge \ell$ crosses $\rho$ and
 every $\tau_{i_j} $ with $j\ge \ell'$ crosses $\rho'$.

\begin{figure}
\input{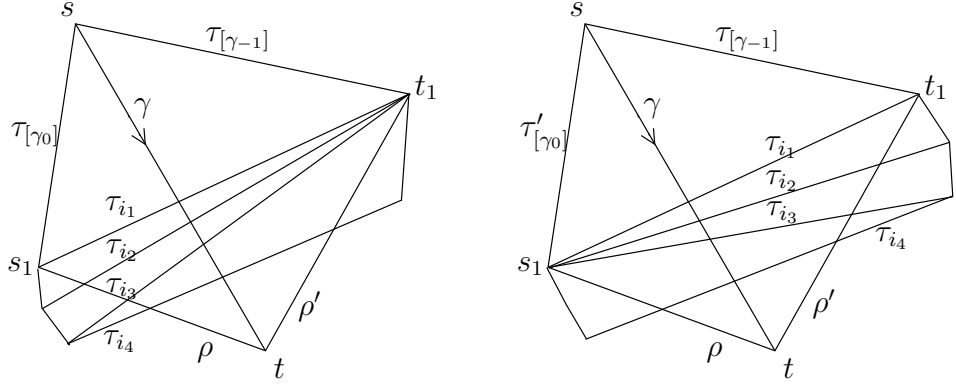}
\caption{Two examples of triangulations. On the left $\ell=2,\ell'=4$
  and on the right $\ell=4,\ell'=2$}\label{fig ell}
\end{figure}

\begin{lem}\label{lem 3}
Let $S$ be simply connected. Let
$T'=\mu_{i_2}\mu_{i_3}\ldots\mu_{i_d}(T)$ and let
$\tilde B_{T'}=(b'_{jk})=\mu_{i_2}\mu_{i_3}\ldots\mu_{i_d}(\tilde B_T)$ be the
  $2n\times n$ matrix obtained by applying the sequence of mutations $\mu_{i_2},\mu_{i_3},\ldots,\mu_{i_d}$ to $\tilde B_T$. Then
\[b'_{(n+i_j)i_1}= \left\{
\begin{array}{ll}
1&\textup{if $j<\ell'$}\\
0&\textup{otherwise.}
\end{array}  \right.\]
\end{lem}

\begin{pf} First note that each $b_{ij}$ is either
  $-1,0$ or $1$, since $S$ is simply connected. Let $Q_T$ be the quiver 
  defined in section \ref{sect quiver}, that is, $Q_T$ has $2n$ vertices labeled
  $1,2,\ldots,2n$, and there is precisely one arrow $i\to j$ whenever
  $\tilde b_{ij}=1$. 
 Consider the following full
  subquiver of $Q_T$
\[ 
\xymatrix@C=2pt@R=20pt{
&&&{[\zg_1]} \ar@{-}[dr]\ar@{-}[dl]
&&[\zg_2] \ar@{-}[dr]\ar@{-}[dl]
&&\ldots&&\ldots 
&&[\zg_{j-1}] \ar@{-}[dr]\ar@{-}[dl]\\
\zG_j &=&i_1 \ar@{-}[rr] &&
i_2 \ar@{-}[rr] &&
i_3 \ar@{-}[rr] &&
\ldots\ar@{-}[rr] && 
i_{j-1}\ar@{-}[rr]  &&
i_j \\
&&n+i_1 \ar[u] &&
n+i_2 \ar[u] &&
n+i_3 \ar[u] &&
\ldots&&
n+i_{j-1} \ar[u] &&
n+i_j \ar[u] 
}
\]
where the orientations of the edges without arrowheads is such that
each triangle $i_k, i_{k+1}, [\zg_k]$ is an oriented $3$-cycle and
$i_{k-1}\ot i_{k}$ for all $k<\ell'$ and $i_{\ell'-1}\to i_{\ell'}$. 
Moreover, the vertices labeled $[\zg_{k}]$ are present only if  $\tau_{[\zg_{k}]}$ is an interior arc of $(S,M)$, thus, whenever  $\tau_{[\zg_{k}]}$ is a boundary arc in $(S,M)$, one has to remove the vertex  $[\zg_{k}]$ and all incident arrows from the above diagram. 
Since $S$ is simply connected, $\zG_j$ is a full subquiver of
$\mu_{j+1}\ldots\mu_d(Q_T)$. Therefore, it suffices to show the following 
\begin{itemize}
\item[Claim:]\emph{ There is precisely one arrow $(n+i_j)\to i_1$  between
    $(n+i_j)$ and $ i_1$ in
    $\zG_j'=\mu_{2}\ldots\mu_{i_j}(\zG_j)$ if  for each
    $k=2,\ldots,j$ the orientation of the 
    arrow between $i_{k-1}$ and $i_k$ is $i_{k-1}\ot i_{k} $ in $\zG_j$;
 and there is no 
    arrow between $(n+i_j)$ and $i_1$ in $\zG_j'$ otherwise.}
\end{itemize}
If the orientations of the arrow between $i_{j-1}$ and $i_j$ is 
$i_{j-1} \to i_j $ in $\zG_j$, then there is no arrow between $(n+i_j)$  and
    $i_{j-1}$ in the quiver $\mu_{i_j}(\zG_j)$, and the subsequent
    mutations $\mu_{i_{j-1}},\ldots,\mu_{i_2}$ will not change the
    arrows at the vertex  $(n+i_j)$. Thus in this case, there are no
    arrows between  $(n+i_j)$ and $i_1$ in $\zG_j'$.
On the other hand, if  the orientation of the arrow between $i_{j-1}$
and $i_j$ is  
 $i_{j-1} \ot i_j $ in $\zG_j$, then there is an arrow  $(n+i_j)\to
    i_{j-1}$ in $\mu_{i_j}(\zG_j)$. Moreover, by deleting the vertices
    $i_j,[\zg_{j-1}]$ and 
    $n+i_{j-1}$, we see that  the quiver $\mu_{i_j}(\zG_j)$ contains the
    following quiver as a full subquiver.
\[ 
\xymatrix@C=2pt@R=20pt{
&&&{[\zg_1]} \ar@{-}[dr]\ar@{-}[dl]
&&[\zg_2] \ar@{-}[dr]\ar@{-}[dl]
&&\ldots&&\ldots 
&&[\zg_{j-2}] \ar@{-}[dr]\ar@{-}[dl]\\
\zG_{j-1} &\cong&i_1 \ar@{-}[rr] &&
i_2 \ar@{-}[rr] &&
i_3 \ar@{-}[rr] &&
\ldots\ar@{-}[rr] && 
i_{j-2}\ar@{-}[rr]  &&
i_{j-1} \\
&&n+i_1 \ar[u] &&
n+i_2 \ar[u] &&
n+i_3 \ar[u] &&
\ldots&&
n+i_{j-2} \ar[u] &&
n+i_j \ar[u] 
}
\]
The claim now follows by induction.
\qed
\end{pf}

We keep the setup of Lemma \ref{lem 3}. 
Define a mapping $f$ for any complete $(T,\zg)$-path $\za=(\za_1,\za_2,\ldots,\za_{2d+1})$
 by 

\[ f(\za) =\left\{
\begin{array}{rcl}
(\za_{2\ell-1},\za_{2\ell},\ldots,\za_{2d+1}) & \textup{if}&\za_1={\tau_{[\zg_{-1}]}}
  \\
(\za_{2\ell'-1},\za_{2\ell'},\ldots,\za_{2d+1}) & \textup{if}&\za_1={\tau_{[\zg_0]}}.
  \\
\end{array}  \right.
\]

\begin{lem}\label{lem 2}
$f$ is a bijection \[\{\textup{complete $(T,\zg)$-paths}\} \stackrel{f}{\to} \{
  \textup{complete $(T,\rho)$-paths}\} \cup  \{\textup{complete
  $(T,\rho')$-paths}\},\]
where $f(\za)$ is a complete $(T,\rho)$-path if
  $\za_1={\tau_{[\zg_{-1}]}}$, and $f(\za) $ is a complete
  $(T,\rho')$-path if $\za_1={\tau_{[\zg_0]}}$. Moreover
\[
\left. \begin{array}{rcl}
{\displaystyle \frac{x_{\tau_{[\zg_{-1}]}}}{x_{{i_1}}}}\, x(f(\za))&=& x(\za)\\ \ \\
y(f(\za))&=& y(\za)
\end{array}  \right\} \qquad \textup{if }  \za_1={\tau_{[\zg_{-1}]}} \]
\ \\
\[\left.\begin{array}{rcl}
{\displaystyle \frac{x_{\tau_{[\zg_0]}}}{x_{{i_1}}}}\, x(f(\za))&=& x(\za)\\ \ \\
Y'\,y(f(\za))&=& y(\za)
\end{array}  \right\}\qquad \textup{if }  \za_1={\tau_{[\zg_0]}}
\]
where $Y'= y_{i_1}y_{i_2}\ldots y_{i_{\ell'-1}}$.
\end{lem}

\begin{pf} First we show that $f$ is well defined.
Clearly, $f(\za)$ is a $T$-path. If $\za_1={\tau_{[\zg_{-1}]}}$, then the even
arcs $\tau_{i_\ell},\tau_{i_{\ell+1}},\ldots,\tau_{i_d}$ of $f(\za)$
are crossing $\rho$ in order,
and, by Lemma \ref{lem 5}, $\rho$ does not cross any other arcs of $T$.
On the other hand, if $\za_1={\tau_{[\zg_0]}}$,
then the even 
arcs $\tau_{i_{\ell'}},\tau_{i_{\ell'+1}},\ldots,\tau_{i_d}$ of $f(\za)$
are crossing $\rho'$ in order, and, again by Lemma \ref{lem 5},
$\rho'$ does not cross any other arcs of $T$. This shows that
$f(\za)$ satisfies 
axiom (T1). The axiom (T2) follows from Lemma \ref{lem 0} (b), and therefore,
 $f$ is well defined.

Next we show that $f$ is injective. Suppose $f(\za)=f(\za')$. Then it
follows from the definition of $f$ that $\za_1=\za_1'$ and 
$\za_2=\za_2'$. Suppose first that $\za_1={\tau_{[\zg_{-1}]}}$, see
Figure \ref{fig ell}. 
If $\ell=2$ then
$f(\za)=(\za_3,\za_4,\ldots,\za_{2d+1})$, while
$f(\za')=(\za'_3,\za'_4,\ldots,\za'_{2d+1})$, and thus
$\za=\za'$. Otherwise, $\ell'=2$ and since the even arcs of $\za$ and $\za'$ must cross $\zg$
in order, there is only one possibility for the 
first $2\ell-2$ arcs in $\za$ and $\za'$, namely
$\za_1=\za_1'=\tau_{[\zg_{-1}]}$,
$\za_2=\za_2'=\tau_{i_1}$ and 
$\za_{2j-1}=\za_{2j}=\tau_{i_j}=\za'_{2j-1}=\za'_{2j}$ 
for $j=2,3,\ldots,\ell-1$, see the right hand side of Figure \ref{fig ell}. Hence $\za=\za'$. Suppose now that
$\za_1={\tau_{[\zg_0]}}$. If $\ell' =2$ then
clearly $\za=\za'$. Otherwise, there is only one possibility for the
first $2\ell'-2$ arcs in $\za$ and $\za'$, namely
 $\za_{2j-1}=\za_{2j}=\tau_{i_j}$ and  $\za'_{2j-1}=\za'_{2j}=\tau_{i_j}$
for $j=2,3,\ldots,\ell'-1$, see the left hand side of Figure \ref{fig ell}. Hence $\za=\za'$.

It remains to show that $f$ is surjective. Let $\hat\za=(\za_{2\ell-1},
\za_{2\ell},\ldots,\za_{2d+1})$ be a complete $(T,\rho)$-path. Then  let
\[
\za=
\left\{\begin{array}{ll}
({\tau_{[\zg_{-1}]}},\tau_{i_1},\za_{2\ell-1},\za_{2\ell},\ldots,\za_{2d+1})&\textup{if
  $\ell=2$}\\ 
({\tau_{[\zg_{-1}]}},\tau_{i_1},\tau_{i_2},\tau_{i_2},\tau_{i_3},\tau_{i_3},\ldots,
\tau_{i_{\ell-1}},\tau_{i_{\ell-1}},
	\za_{2\ell-1},\za_{2\ell},\ldots,\za_{2d+1})&\textup{if
  $\ell>2$}.
\end{array}\right.\]
Then $\za$ is a $T$-path,
 its even arcs are $\tau_{i_1},\tau_{i_2},\ldots,\tau_{i_d}$ in
order, and then it follows from Lemma \ref{lem 0} (b) that $\za$ is a 
complete $(T,\zg)$-path. Moreover, $f(\za)=\hat\za$.
On the other hand, let  $\check\za=(\za_{2\ell'-1},
\za_{2\ell'},\ldots,\za_{2d+1})$ be a complete $(T,\rho')$-path. Then
let
\[\za=
\left\{\begin{array}{ll}
({\tau_{[\zg_0]}},\tau_{i_1},\za_{2\ell'-1},\za_{2\ell'},\ldots,\za_{2d+1})&\textup{if
  $\ell'=2$}\\ 
({\tau_{[\zg_0]}},\tau_{i_1},\tau_{i_2},\tau_{i_2},\tau_{i_3},\tau_{i_3},\ldots,
\tau_{i_{\ell'-1}},\tau_{i_{\ell'-1}},
        \za_{2\ell'-1},\za_{2\ell'},\ldots,\za_{2d+1})&\textup{if
  $\ell'>2$}.
\end{array}\right.\]
Then $\za$ is a
complete $(T,\zg)$-path, and $f(\za)=\check\za$. Thus $f$ is surjective, and hence a bijection.

Now we show the identities in the statement of Lemma \ref{lem 2}.
Suppose first that $\za_1={\tau_{[\zg_{-1}]}}$. Then $\za_2=\tau_{i_1}$ is not $\zg$-oriented. If $\ell=2$
then $f(\za)=(\za_3,\za_4,\ldots,\za_{2d+1})$ and
\[
 x(\za) = {\displaystyle \frac{x_{\tau_{[\zg_{-1}]}}}{x_{{i_1}}}}\ x({\scriptstyle
 f(\za)})\qquad 
 \textup{and} \qquad  y(\za)=  y({\scriptstyle f(\za)}).
\]
On the other hand, if $\ell\ne 2$ then 
\[\za=({\tau_{[\zg_{-1}]}},\tau_{i_1},\tau_{i_2},\tau_{i_2},\ldots,\tau_{i_{\ell-1}},
\tau_{i_{\ell-1}},\za_{2\ell-1},\za_{2\ell},\ldots)\]
where none of the first $\ell-1$ even arcs is $\zg$-oriented, see the right hand side of Figure \ref{fig ell},
and $f(\za)=(\za_{2\ell-1},\za_{2\ell},\ldots).$ Hence
\[
 x(\za)= {\displaystyle\frac{x_{\tau_{[\zg_{-1}]}}}{x_{{i_1}}}}\ x({\scriptstyle f(\za)})\qquad
 \textup{and} \qquad y(\za)= y({\scriptstyle f(\za)}). 
\]
This proves the lemma in the case $\za_1=\tau_{[\zg_{-1}]}$.

Suppose now that $\za_1={\tau_{[\zg_0]}}$. Then $\za_2$ is $\zg$-oriented. If $\ell'=2$
then $f(\za)=(\za_3,\za_4,\ldots,\za_{2d+1})$ and 

\[
 x(\za) ={\displaystyle\frac{x_{{\tau_{[\zg_0]}}}}{x_{{i_1}}}} \ x({\scriptstyle f(\za)})\qquad
 \textup{and} \qquad
 y(\za)= y_{i_1}\ y({\scriptstyle f(\za)}) = Y'\ y({\scriptstyle f(\za)}).
\]
On the other hand, if $\ell'\ne 2$ then 
\[\za=({\tau_{[\zg_0]}},\tau_{i_1},\tau_{i_2},\tau_{i_2},\ldots,\tau_{i_{\ell'-1}},
\tau_{i_{\ell'-1}},\za_{2{\ell'-1}},\za_{2{\ell'}},\ldots)\]
where each of the first $\ell'-1$ even arcs is $\zg$-oriented, see the left hand side of Figure \ref{fig ell},
and $f(\za)=(\za_{2{\ell'-1}},\za_{2{\ell'}},\ldots).$ Hence
\[
 x(\za) = \frac{x_{{\tau_{[\zg_0]}}}}{x_{{i_1}}}\ x({\scriptstyle f(\za)})\qquad \textup{and}
 \qquad 
 y(\za)=y_{{i_1}} y_{{i_2}} \ldots y_{{i_{\ell'-1}}}\ y({\scriptstyle f(\za)})
 =Y'\ y\left({\scriptstyle f(\za)}\right).
\]
\qed
\end{pf}

\paragraph{\bf Proof of Theorem \ref{thm main} in the simply
    connected case.}
We prove the Theorem by induction on $d$. If $d=0$ then $\zg\in T$ and the only complete $(T,\zg)$-path is the path $\za=(\zg)$. Hence $\sum_\za x(\za)y(\za)=x_\zg$, as desired. Suppose now $d\ge 1$.
By Corollary \ref{cor 1}, the triangulation
$T'=\mu_{i_2}\mu_{i_3}\ldots\mu_{i_d}(T)$ contains
${\tau_{[\zg_{-1}]}},{\tau_{[\zg_0]}},\rho,\rho',$ and $\tau_{i_1}$,
and we thus have the exchange relation
\begin{equation}\label{eq 96}
x_{i_1}x_\zg = Y\, x_{\tau_{[\zg_{-1}]}} x_\rho+ Y'\,x_{{\tau_{[\zg_0]}}}x_{\rho'}
\end{equation}
where 
\[Y=\prod_{b'_{(n+i_j)i_1}\le 0} y_{i_j}^{-b'_{(n+i_j)i_1}} \qquad \textup{and} \qquad
 Y'=\prod_{b'_{(n+i_j)i_1}\ge 0} y_{i_j}^{b'_{(n+i_j)i_1}}\]
 and where $B'=(b'_{ij})$ is
 the $2n\times n$ matrix given by $B'=\mu_{i_2}\ldots\mu_{i_d}(B_T).$
By Lemma \ref{lem 3}, we have $b'_{(n+i_j)i_1}=1$ if 
$j<\ell'$  and $b'_{(n+i_j)i_1}=0$ otherwise. In
 particular, $Y=1$; and $Y'=y_{i_1}y_{i_2}\cdots y_{i_{\ell'-1}}$ as in
 Lemma \ref{lem 2}.
 
The number of crossings between $\zg$ and $T$ is strictly larger than
the number of crossings between $\rho $ and $T$  and strictly larger than
the number of crossings between $\rho' $ and $T$. 
Hence, by induction, equation (\ref{eq 96}) implies

\begin{equation}\label{eq*}
 x_\zg = \sum_\zb \frac{x_{\tau_{[\zg_{-1}]}}}{x_{i_1}}\, x(\zb) \, y(\zb)+
 \sum_{\zb'} \frac{x_{{\tau_{[\zg_0]}}}}{x_{i_1}} \, x(\zb')\, Y'\, y(\zb'),
\end{equation}
where, in  the first sum, $\zb$ runs over all complete $(T,\rho)$-paths
in $(S,M)$, and, in the second sum, $\zb'$ runs over all complete
$(T,\rho')$-paths in $(S,M)$.

Applying  Lemma \ref{lem 2} to equation (\ref{eq*})
we get 
\begin{equation}\label{eq**}
x_\zg = \sum_{\za:\za_1={\tau_{[\zg_{-1}]}}} x(\za)y(\za) +\sum_{\za : \za_1 ={\tau_{[\zg_0]}}} x(\za)y(\za),
\end{equation}  
where, in the first sum, $\za$ runs over all complete $(T,\zg)$-paths
starting with ${\tau_{[\zg_{-1}]}}$, and, in the second sum,  $\za$ runs over all
$(T,\zg)$-paths starting with ${\tau_{[\zg_0]}}$.

Since $S$ is simply connected and $\za_2=\tau_{i_1}$ for any complete
$(T,\zg)$-path $\za$, Lemma \ref{lem 0} (b) implies that we must have
$\za_1={\tau_{[\zg_{-1}]}}$ or $\za_1 = {\tau_{[\zg_0]}}$ for any complete $(T,\zg)$-path $\za$.

 Therefore the right hand side of
equation (\ref{eq**}) equals the right hand side of the formula in
Theorem \ref{thm main}. This proves  the
  theorem in the case where $S$ is simply connected.
The general case will be proved in section \ref{sect general case}.

\end{section}


\begin{section}{Universal cover}\label{sect cover}

In this section, we will use covering techniques to prove Theorem \ref{thm main} in the general case. The idea of the proof is to work in the universal cover of $S$ in which we can apply the Theorem \ref{thm main} for simply connected surfaces, and then project the results to the original surface $S$. First we need to investigate Galois coverings of quivers.

\begin{subsection}{Galois coverings}

Let $k$ be a field.
Any quiver $Q$ can be considered as a $k$-category, whose objects are
the vertices of $Q$ and for two vertices $x,y$ of $Q$, the space
$\Hom_Q(x,y)$ is spanned by the set of arrows from $x$ to $y$ in $Q$.
The composition of morphisms is the composition of arrows.

Let $Q,\tilde Q$ be two quivers and 
$\pi:\tilde Q \to Q$ a functor, that is, $\pi $ maps vertices to
vertices and arrows to arrows respecting the composition.
Let $\textup{Aut}(\tilde Q)$ be the group of autoequivalences of
$\tilde Q$, and let $G$ be a subgroup of $\textup{Aut}(\tilde Q)$ whose
action on the objects of $\tilde Q$ is free  (that is $g(x)\ne x$ for
all vertices $x$ and all $g\in G$, $g\ne id$) 
and such that $\pi\circ g=\pi$ for all $g\in G$.

Then $\pi:\tilde Q\to Q$ is called a \emph{Galois covering with group
  $G$} if  the orbit category $\tilde Q/G$ of $\tilde Q$ under
  the action of $G$ is equivalent to $Q$.

If $j$ is a vertex of a quiver $Q$, let $\mu_j$ denote the mutation in
direction $j$. Let $\pi:\tilde Q\to Q$ be a Galois covering of quivers
with group $G$. The elements of the fiber $\pi^{-1}(j)$ of $j$ in
$\tilde Q$ are of the form $j_g$ with $g\in G$. Let $\tilde \mu_j$
denote the infinite sequence of mutations $(\mu_{j_g})_{g\in G}$.

\begin{lem}\label{lem galois}
If $\pi:\tilde Q\to Q$ is a Galois covering of quivers with group
$G$, then there exists a Galois covering 
$\pi':\tilde\mu_j(\tilde Q)\to \mu_j(Q)$  of
the mutated quivers with the same group $G$.
\end{lem}

\begin{pf}
The mutated quiver $\mu_j(Q)$ is obtained from $Q$ by the following
operations: 
\begin{enumerate}
\item For every path $i\stackrel{\za}{\to} j\stackrel{\zb}{\to} k$
  introduce a new arrow $i\stackrel{[\za\zb]}{\longrightarrow} k$. 
\item Replace each arrow $\za$ starting at (respectively ending at)
  $j$ by its opposite $\za^*$ ending at (respectively starting at) $j$.
\item Cancel all $2$-cycles.
\end{enumerate}
For each arrow $i\stackrel{\za}{\to} j$ (respectively
$j\stackrel{\zb}{\to} k$) and for each $g\in G$ there exists a unique arrow 
$i_h\stackrel{\za_g}{\to} j_g$ (respectively
$j_g\stackrel{\zb_g}{\to} k_{h'}$) that ends at $j_g$ (respectively
starts at $j_g$) for some $h,h'\in G$, and such that $\pi(\za_g)=\za$.

Thus the mutated quiver $\tilde\mu_j(\tilde Q)$ is obtained from $\tilde Q$
by the following 
operations:
 \begin{enumerate}
\item For every path $i_h\stackrel{\za_g}{\to} j_g\stackrel{\zb_g}{\to} k_{h'}$
with $h,g,h'\in G$  introduce a 
new arrow $i_h\stackrel{[\za\zb]_g}{\longrightarrow} k_{h'}$.
\item Replace each arrow $\za_g$ starting at (respectively ending at)
  $j_g$ by  its opposite $\za^*_g$ ending at (respectively starting at)
  $j_g$.
\item Cancel all $2$-cycles.
\end{enumerate}   
Define $\pi':\tilde\mu_j(\tilde Q)\to\mu_j(Q)$ as follows. 
Let  $\pi'=\pi$ on all vertices
of $\tilde\mu_j \tilde Q$ and let $\pi'=\pi$ on those arrows that lie in both
$\tilde\mu_j(\tilde Q)$ and $\tilde Q$; and let $\pi'([\za\zb]_g)=[\za\zb]$
and $\pi'(\za^*_g)=\za^*$. The group $G$ acts on the vertices of $\tilde \mu_j(\tilde Q)$
in the same way as on $\tilde Q$. Its action on the arrows of $\tilde \mu_j(\tilde Q)$ that have not been modified by the mutation $\tilde \mu_j$ is the same as on $\tilde Q$. On the new arrows $[\za\zb]_g$, the action of an element $g'\in G$ is given by 
\[ g'([\za\zb]_g) = \left(i_{g'h}\stackrel{[\za\zb]_{g'g}}{\longrightarrow}k_{g'h}\right), 
\]
for all $g,g',h,h'\in G$. Clearly, $\pi'\circ g=\pi'$ for all $g\in G$, and thus $\pi'$ is a Galois covering of
quivers with the same group as $\pi$.
\qed
\end{pf}  

\end{subsection} 

\begin{subsection}{Triangulations and quivers}\label{sect TQ}

Let $\pi:\tilde S \to S$ be a universal cover of the surface $S$, and
let $\tilde M=\pi^{-1}(M)$ and $\tilde T=\pi^{-1}(T)$.

From the theory of covering spaces, we know that for each point $p\in
S$ the fiber $\pi^{-1}(p)$ is the fundamental group $\Pi_1(S)$ of the
surface $S$. Hence $\pi^{-1}(p)=\{\tilde p_\zs\mid\zs\in\Pi_1(S)\}$.
Furthermore, for each path $\tau $ in $S$ with starting point $p\in
S$ and for each $\zs\in \Pi_1(S)$
there is a unique lift $\tilde \tau_\zs$ in $\tilde S$ with starting 
point $\tilde p_\zs$ in the fiber of $p$. 
The triangulation $\tilde T$ consists of arcs $\tilde \tau_\zs$ where
$\zs\in\Pi_1(S)$ and $\pi(\tilde\tau_\zs)\in T$.

Let $Q_T$ be the quiver associated to the triangulation $T$ as in section \ref{sect quiver}. Recall
that $Q_T$ has $2n$ vertices labeled $1,2,\ldots,2n$, where the first
$n$ vertices correspond to the interior arcs of the triangulation $T$,
and the vertex $n+i$ correspond to the coefficient $y_i$, for
$i=1,2,\ldots,n$. Since the cluster algebra  $\mathcal{A}$ has principal
coefficients in the initial seed, we have an arrow $(n+i)\to j$ in $Q_T$
 if $i=j$ and there are no arrows between $(n+i)$ and $j$ in $Q_T$ if
$i\ne j$.

Let $Q_{\tilde T}$ be the (infinite) quiver associated in the same way to
triangulation $\tilde T$ of $(\tilde S,\tilde M)$. 

The vertices of the quiver $Q_{\tilde T}$ are labeled by tuples
$(\tau,\zs)$ where $\tau\in \{1,2,\ldots,2n\}$ is a vertex of $Q_T$ and
$\zs\in\Pi_1(S)$, and the arrows by $\za_\zs$ where $\za$ is an arrow 
in $Q_T$ and $\zs\in \Pi_1(S)$.
Hence $\pi$ induces a functor $\pi: Q_{\tilde T}\to Q_T$. The fundamental
group $ \Pi_1(S)$ acts on $Q_{\tilde T}$ by autoequivalences as 
$\rho\cdot (\tau,\zs)=(\tau,\rho\zs)$, for all $\tau\in\{1,2,\ldots,2n\} $, and all $\zs,\rho\in
  \Pi_1(S)$; and $\rho\cdot\za_\zs =\za_{\rho \zs}$, for all arrows $\za$ in 
  $Q_T$, and all $\zs,\rho\in  \Pi_1(S)$.
This action is free on vertices and for all $\zs\in\Pi_1(S)$, we have
$\pi\circ\zs=\pi$. Moreover, $Q_T\cong Q_{\tilde T} /\Pi_1(S)$, and thus
 $\pi:Q_{\tilde T}\to Q_T$ is a Galois covering of quivers
with group $\Pi_1(S)$. Applying Lemma \ref{lem galois}, we get the
following theorem.

\begin{thm}\label{thm galois}
Let $T'$ be the triangulation of $(S,M)$ obtained from $T$ by a
sequence of mutations $T'=\mu_{j_1}\mu_{j_2}\ldots \mu_{j_s} (T)$. 
Let $\tilde T=\pi^{-1}(T)$ be the lifted triangulation of $(\tilde
S,\tilde M)$ and denote by $Q_{\tilde T}$ the corresponding quiver. 
Let $\tilde
T'=\tilde\mu_{j_1}\tilde\mu_{j_2}\ldots\tilde\mu_{j_s}(\tilde T)$ 
 and $Q_{\tilde  T'}
 =\tilde\mu_{j_1}\tilde\mu_{j_2}\ldots\tilde\mu_{j_s}( Q_{\tilde T})$.
Then $T'=\pi(\tilde T')$ and $ \pi : Q_{\tilde T'}\to Q_{T'}$ is a Galois covering of quivers with group $\Pi_1(S)$.
\end{thm}  
\qed
\end{subsection}

\begin{subsection}{Cluster algebras}\label{sect CA}
In this subsection, we will use the universal cover to define a new cluster algebra $\tilde{\mathcal{A}}$ in which we can use Theorem \ref{thm main} for simply connected surfaces. We will not work inside the universal covering space $\tilde S$ itself but we will restrict to a subsurface of $S$ that contains a lift of the arc $\zg$ and has a finite triangulation. We also need a special coefficient system that takes care of the principal coefficients as well as of the boundary of the surface.

We keep the setup of section \ref{sect expansion}.
Let $\pi:\tilde S \to S$ be a universal cover of the surface $S$, and
let $\tilde M=\pi^{-1}(M)$ and $\tilde T=\pi^{-1}(T)$. 

Choose a point $\tilde s $ in the fiber $\pi^{-1}(s)$ of the starting point $s$ of the arc $\zg$. There exists a unique lift $\tilde
\zg$ of $\zg$ starting at $\tilde s$. Then $\tilde \zg$ is the
concatenation of subpaths $\tilde \zg_0,\tilde \zg_1,\ldots,\tilde
\zg_{d}$ where $\tilde \zg_k$ is a path from a point $\tilde p_k$ to
a point $\tilde p_{k+1}$ such that $\tilde \zg_k$ is a lift of $\zg_k$
and $\tilde p_k\in\pi^{-1}(p_k)$, for $k=0,1,\ldots, d$.
Let $\tilde t=\tilde p_{d+1}\in \pi^{-1}(t)$.

For $k$ from $1$ to $d$, let $\tilde \tau_{i_k}$ be the unique lift of
$\tau_{i_k}$ running through $\tilde p_k$. For $k$ from $1$ to
$d-1$ let $\tilde 
\tau_{[\zg_k]}$ be the  unique lift of $\tau_{[\zg_k]}$ that is
bounding a triangle $\tilde \zD_k$ in $\tilde T$ with $\tilde
\tau_{i_k}$ and $\tilde \tau_{i_{k+1}}$. Let  $\tilde
\tau_{[\zg_{-1}]} $ and $\tilde\tau_{[\zg_{0}]}$ be the unique lifts of
    ${\tau_{[\zg_{-1}]}}$ and $\tau_{[\zg_{0}]}$ that, together with $\tilde\tau_{i_1}$, are bounding a triangle $\tilde\zD_0$,  and let  $\tilde
\tau_{[\zg_{d}]}$ and $\tilde\tau_{[\zg_{d+1}]}$ be the unique lifts of
    $\tau_{[\zg_{d}]}$ and $\tau_{[\zg_{d+1}]}$ that, together with $\tilde\tau_{i_d}$, are bounding a
triangle $\tilde\zD_{d}$.

\begin{definition}
Let $\tilde S(\zg)\subset \tilde S$ be the union of the $d$ triangles
$\tilde \zD_0, \tilde \zD_1,\ldots,\tilde
\zD_{d}$ and let $\tilde
M(\zg)=\tilde M\cap \tilde S(\zg)$ and $\tilde T(\zg)=\tilde T\cap
\tilde S(\zg)$.
\end{definition}

\begin{prop}\label{prop 1} 
$(\tilde S(\zg),\tilde M(\zg))$ is a simply connected unpunctured
surface of which $\tilde T(\zg) $
is a triangulation. This 
triangulation $\tilde T(\zg)$ consists of $d$ interior arcs 
 and $d+3$ boundary arcs.

In particular,  each triangle $\tilde \zD_k$ in $\tilde T(\zg)$
contains a boundary arc.
\end{prop}  
\begin{pf} This follows immediately from the construction.
\qed
\end{pf}  

Associate a quiver $\tilde Q_{\tilde T(\zg)}$ to $\tilde T(\zg)$ as follows:
The set of vertices of $\tilde Q_{\tilde T(\zg)}$ is
\[\{
i_1,i_2,\ldots,i_d,c_1,c_2,\ldots,c_d,
[\zg_{-1}],[\zg_{0}],\ldots,[\zg_{d+1}]
\}\]
corresponding to the $d$ interior arcs
$\tilde \tau_{i_1},\tilde \tau_{i_2},\ldots,\tilde \tau_{i_d}$ of  $\tilde T(\zg)$, 
the $d+3$ boundary arcs $\tilde \tau_{[\zg_{-1}]},\tilde \tau_{[\zg_0]},\ldots, \tilde \tau_{[\zg_{d+1}]}$ of  $\tilde
T(\zg)$, and one extra vertex $c_j$ for each 
interior arc $\tilde \tau_{i_j}$.
The set of arrows consists of one arrow $ {c_j}\to {i_j}$, for each $j$, $1\le j \le d$; and one arrow $i_j\to i_k$ (respectively   $i_j\to [\zg_\ell]$; $ [\zg_\ell]\to i_j$ ) whenever $(\tilde\tau_{i_j},\tilde\tau_{i_k})$ (respectively $(\tilde\tau_{i_j},\tilde\tau_{[\zg_\ell]})$) are sides of the same triangle in $\tilde T(\zg)$ with $\tilde \tau_{i_k}$ following $\tilde \tau_{i_j}$ (respectively 
$\tilde \tau_{[\zg_\ell]}$ following $\tilde \tau_{i_j}$;  $\tilde \tau_{i_j}$ following $\tilde \tau_{[\zg_\ell]}$) in the counter-clockwise order.
Let us point out that the difference between this construction and the one in section \ref{sect quiver} is that it involves also the boundary arcs of the surface.

Since each of the triangles in $\tilde S(\zg)$ is a lift of an oriented
triangle in $S$, the orientation of $S$ lifts to an orientation of
$\tilde S(\zg)$ and we have 

\begin{lem}\label{lem or1} A path $\tilde \za_{2k}$ along
  the arc $\tilde\tau_{i_k}$ is $\tilde {\zg}$-oriented if and only if
  $\pi(\tilde{\za}_{2k})$ along the arc $\tau_{i_k}$ is $\zg$-oriented.
\qed\end{lem}

Let 
\[\tilde{\mathbb{P}}=\textup{Trop}(u_{i_1},u_{i_2},\ldots,u_{i_d},
u_{[\zg_{-1}]}, u_{[\zg_{0}]},\ldots , u_{[\zg_{d+1}]} )
\]
be the tropical semifield on $2d+3$ generators, where the first $d$
generators correspond to the $d$ interior arcs $\tilde\tau_{i_k},
k=1,2,\ldots,d$ and the  last $d+3 $ generators $ u_{[\zg_{j}]}$
correspond to the $d+3$ boundary arcs $\tilde \tau_{[\zg_j]},
j=-1,0,\ldots,d+1$.
Let $\tilde{\mathcal{F}}$ be the field of rational functions in $d$ variables with coefficients in $\mathbb{Q}\tilde{\mathbb{P}}$.

\begin{definition} \label{def cl} Let
$\tilde{\mathcal{A}}=\mathcal{A}(\tilde{\mathbf{x}},\tilde{\mathbf{y}},\tilde
B)$ be the cluster algebra in $\tilde{\mathcal{F}}$ given by the initial
seed $\tilde \Sigma$:
\[
\begin{array}{rl}
\tilde{\mathbf{x}} = \{\tilde{x}_{i_1},\ldots,\tilde{x}_{i_d}\}  
& \textup{where } \tilde x_{i_k}=x_{\tilde \tau_{i_k}} \\
\tilde{\mathbf{y}} = \{\tilde{y}_{i_1},\ldots,\tilde{y}_{i_d}\}
& \textup{where } \tilde y_{i_k}=u_{i_k} \prod_{[\zg_j]\to i_k \in
  \tilde Q_{\tilde T(\zg)}} u_{[\zg_j]}
\prod_{[\zg_j]\ot i_k \in
  \tilde Q_{\tilde T(\zg)}} u_{[\zg_j]}^{-1}\\
\tilde B=B_{\tilde T(\zg)}
\end{array}  
\]
\end{definition}

Thus the cluster algebra has coefficients combined from principal
coefficients, the $u_{i_k}$ term, and
boundary coefficients, the $u_{[\zg_j]}$ terms coming from boundary arcs.

For each $k=1,2,\ldots,d$, denote by $\tilde x_{i_k}'$ the cluster variable
obtained by mutation in direction $k$. That is
$\mu_k (\tilde{\mathbf{x}})=\tilde{\mathbf{x}}\setminus\{\tilde x_{i_k}\}
    \cup \{\tilde x'_{i_k}\}$.

\begin{prop}\label{prop 2}
The cluster algebra $\tilde{\mathcal{A}}$ is an acyclic cluster algebra
of finite type $A_{d}$ with acyclic seed $\tilde\Sigma$. In particular,
$\tilde{\mathcal{A}}$ is generated over $\mathbb{Z}\tilde{\mathbb{P}}$ by the $2d$ cluster variables
\[\tilde x_{i_1},\tilde x_{i_2},\ldots,\tilde x_{i_d},\tilde
x'_{i_1},\tilde x'_{i_2},\ldots,\tilde x'_{i_d}.\]
\end{prop}  

\begin{pf}  The surface $\tilde S(\zg)$ is topologically a polygon with
  $d+3$ vertices,  which implies that the cluster algebra   is of type
  $A_{d}$. 
The seed $\tilde \Sigma $ is acyclic since,  by Proposition \ref{prop 1}, each triangle in $\tilde T(\zg)$ has at least one side given by a boundary arc.
 The last statement now
  follows from  \cite[Cor 1.21]{BFZ}.
\qed
\end{pf}

The following Lemma is shown in \cite{ST}. It will allow us to
use induction later.
\begin{lem}\label{lem exchange}\cite[Lemma 4.6]{ST}
Let $T$  be a triangulation of an unpunctured
  surface $(S,M)$, and let
  ${\zb}$ be an arc in $S$ which is not in $T$.
 Let $k$ be the number of crossings
  between ${\zb}$ and $T$. 
Then there exist five arcs
  $\rho_1,\,\rho_2,\,\zs_1,\,\zs_2$ and ${\zb'}$ in $S$ such that 
\begin{itemize}
\item[(a)] each of $\rho_1,\,\rho_2,\,\zs_1, \,\zs_2$  and ${\zb'}$ crosses $T$
  less than $k$ times, 
\item[(b)] $\rho_1,\,\rho_2,\,\zs_1,\,\zs_2$ are the sides of a simply
  connected 
  quadrilateral $V$ in which ${\zb}$ and ${\zb'}$ are the diagonals.
\end{itemize}   
\end{lem}

\begin{thm}\label{thm cover}
The universal cover $\pi:\tilde S \to S$ induces a homomorphism of
algebras $\pi_*:\tilde {\mathcal{A}} \to \mathcal{A}$ defined by
$\pi_*(\tilde x_{i_k})=x_{i_k}$, $\pi_*(\tilde x'_{i_k})=x'_{i_k}$,
$\pi_*(u_{i_k})=y_{i_k}$ and $\pi_*(u_{[\zg_j]})=x_{[\zg_j]}$.
Moreover,
  if $\tilde \zb$ is an interior arc in $(\tilde S(\zg),\tilde M(\zg))$ which is a lift of an arc $\zb$ in $(S,M)$,   then $\pi_*(\tilde x_{\tilde
  \zb})=x_\zb$.
\end{thm}  

\begin{pf} By Proposition \ref{prop 2}, $\tilde{\mathcal{A}}$ is
  generated over $\mathbb{Z}\tilde{\mathbb{P}} $ by the elements 
$\tilde x_{i_k}$, $\tilde x'_{i_k}$, and the  coefficient semifield is generated by the elements $u_{i_k},u_{[\zg_j]}$, where $1\le k\le d$ and $ -1\le j\le d+1$.
Define $\pi_*$ on the generators $\tilde x_{i_k}$, $\tilde x'_{i_k},
u_{i_k},u_{[\zg_j]}$ as in the statement of the theorem, and extend
it to arbitrary 
elements by the homomorphism property. 

To show that $\pi_*$ is well defined, we have to check that it
preserves the relations between the generators. In
$\tilde{\mathcal{A}}$ these relations are the exchange relations (\ref{exchange relation})
\begin{equation}\label{eq 101} \tilde x_{i_k}  \tilde x'_{i_k} =
\left(  \tilde y_{i_k} \prod_{i_j\to
  i_k\in \tilde Q_{\tilde T(\zg)}}  \tilde x_{i_j}
+ \prod_{i_j\ot
  i_k\in \tilde Q_{\tilde T(\zg)}}  \tilde x_{i_j}\right)/(1\oplus
  \tilde y_{i_k}),
\end{equation}
where \[\begin{array}{rcl }
\tilde y_{i_k}&=&u_{i_k} \prod_{[\zg_j]\to i_k \in
  \tilde Q_{\tilde T(\zg)}} u_{[\zg_j]}
\prod_{[\zg_j]\ot i_k \in
  \tilde Q_{\tilde T(\zg)}} u_{[\zg_j]}^{-1},  \quad
  \textup{by equation (\ref{eq 20}), and} \\ 
  \\
  1\oplus
  \tilde y_{i_k} &=& \prod_{[\zg_j]\ot i_k \in
  \tilde Q_{\tilde T(\zg)}} u_{[\zg_j]}^{-1}, \quad  \textup{by definition of $\oplus$} .
\end{array}\]
By definition, $\pi_*$ maps  equation  (\ref{eq 101}) to 
\[ x_{i_k}   x'_{i_k} =  y_{i_k} \prod_{[\zg_j]\to
  i_k\in  Q_{T}}   x_{[\zg_j]} 
\prod_{i_j\to
  i_k\in  Q_{T}}   x_{i_j}
+ \prod_{[\zg_j]\ot
  i_k\in  Q_{T}}   x_{[\zg_j]} 
 \prod_{i_j\ot
  i_k\in Q_{ T}} x_{i_j},
\]
which, using equation (\ref{geometric exchange}) and the fact that  the initial seed in $\mathcal{A}$ has principal coefficients, is easily seen to be an exchange relation in $\mathcal{A}$.

It remains to show that $\pi_*(\tilde x_{\tilde \zb})=x_\zb$ whenever
$\tilde \zb$ is an interior arc that is a lift of $\zb$.
We prove this by induction on  the minimal number  of
crossing points between $\zb $ and $T$.  
If this number is zero, then $\zb\in T$, and $\pi_*(\tilde x_{\tilde\zb}) =x_\zb$, by
definition. 
Otherwise, let $\rho_1,\ \rho_2,\ \zs_1,\ \zs_2,$
and $\zb'$ be as in Lemma \ref{lem exchange}. Suppose  without loss of generality that the relative position of these arcs is as in Figure \ref{fig beta}. Then, in
$\mathcal{A}$, we have the
exchange relation 
\begin{equation}\label{exchange}
 x_{\zb} =(  x_{\rho_1} x_{\rho_2}\,y_\zb^+ + x_{\zs_1}x_{ \zs_2}\,y_\zb^-)/ x_{\zb'},
\end{equation}   
for some coefficients $y_\zb^+$ and $y_\zb^-$.

\begin{figure}
\begin{center}
\input{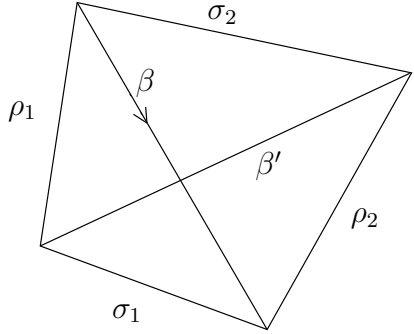}
\caption{The quadrilateral $V$}
\label{fig beta}
\end{center}
\end{figure}

Let $\tilde \rho_1$
and $\tilde \zs_2$ be the unique lifts  of $\rho_1$ and $\zs_2$,
respectively,  that start at the same point as $\tilde \zb$. 
Let $\tilde \zs_1$ and
$\tilde {\zb'}$ be the unique lifts of $\zs_1$ and $\zb'$, respectively, that start at
the endpoint of $\tilde \rho_1$, and let $\tilde \rho_2$ be the unique
lift of $\rho_2$ that starts at the endpoint of $\tilde\zs_2$.
Then $\tilde \rho_1,\,\tilde \rho_2, \,\tilde \zs_1, \,\tilde \zs_2$
form a quadrilateral in $\tilde S(\zg)$ in which $\tilde \zb$ and $\tilde {\zb'}$
are the diagonals. Consequently, in the cluster algebra
$\tilde{\mathcal{A}}$, we have the exchange relation
\begin{equation}\nonumber
 \tilde x_{\tilde \zb} =(  \tilde x_{\tilde \rho_1} \tilde x_{\tilde
  \rho_2} \tilde y_{\tilde \zb}^+ + \tilde x_{\tilde
  \zs_1}\tilde x_{\tilde \zs_2}\tilde y_{\tilde\zb}^-)/ \tilde x_{\tilde {\zb'}}.
\end{equation}   
Therefore 
\begin{equation}\label{eq 95}
 \pi_* (\tilde x_{\tilde \zb} )=
( \pi_*  (\tilde x_{\tilde \rho_1}) \pi_*( \tilde x_{\tilde
   \rho_2})\pi_*(\tilde y_{\tilde \zb}^+)
+  \pi_*(\tilde x_{\tilde  \zs_1})  \pi_*(\tilde x_{\tilde
  \zs_2})\pi_*(\tilde y_{\tilde\zb}^-))/  \pi_*(\tilde x_{\tilde {\zb'}}).
\end{equation}   
In this formula $\tilde y_{\tilde\zb}^+$ is a  product $\tilde y_{\tilde\zb}^+=\tilde u_{\tilde\zb}^+\left( \prod u_{[\zg_j]} \right)$, where  $\tilde u_{\tilde\zb}^+$ is a product of $\tilde u_{i_k}$\!'s, whereas $\prod u_{[\zg_j]}=u_{\tilde\rho_1}^{\zd_1}u_{\tilde\rho_2}^{\zd_2}$
with \[ \zd_i=\left\{\begin{array}{ll}
      0&\textup{if $\tilde \rho_i$ is an interior arc in $(\tilde S(\zg),\tilde M(\zg))$,}    \\
      1& \textup{if $\tilde \rho_i$ is a boundary arc in $(\tilde S(\zg),\tilde M(\zg))$.}
\end{array}\right.\]
Note that if $\zd_i=1$ then $\tilde x_{\tilde \rho_i}=1$.
Similarly,  $\tilde y_{\tilde\zb}^-=\tilde u_{\tilde\zb}^-\left( \prod u_{[\zg_j]} \right)$, where  $\tilde u_{\tilde\zb}^-$ is a product of $\tilde u_{i_k}$\!'s, whereas $\prod u_{[\zg_j]}=u_{\tilde\sigma_1}^{\ze_1}u_{\tilde\sigma_2}^{\ze_2}$
with \[ \ze_i=\left\{\begin{array}{ll}
      0&\textup{if $\tilde \zs_i$ is an interior arc in $(\tilde S(\zg),\tilde M(\zg))$, }   \\
      1& \textup{if $\tilde \zs_i$ is a boundary arc in $(\tilde S(\zg),\tilde M(\zg))$.}
\end{array}\right.\]
Note that if $\ze_i=1$ then $\tilde x_{\tilde \zs_i}=1$.
By induction, equation (\ref{eq 95}) yields
\begin{equation}\nonumber
 \pi_* (\tilde x_{\tilde \zb} )=
( x_{\rho_1} x_{\rho_2} \pi_*(\tilde u_{\tilde \zb}^+)
+x_{\zs_1}  x_{  \zs_2}  \pi_*(\tilde u_{\tilde {\zb}}^-))/ x_{{\zb'}},
\end{equation}   
which, by equation (\ref{exchange}), is equal to $x_{\zb}$ if and only
if 
\begin{equation}\label{claim}
  \pi_*(\tilde u_{\tilde \zb}^+) = y_{\zb}^+
\quad \textup{and}\quad
  \pi_*(\tilde u_{\tilde \zb}^-) = y_{\zb}^-.
\end{equation}

In order to show equation (\ref{claim}), it does not suffice to work with
the triangulations alone, but one has to consider the corresponding
quivers. We will also need to work in the universal cover $(\tilde
S,\tilde M)$ and its triangulation $\tilde T=\pi^{-1} (T)$ instead of the surface $(\tilde S(\zg),\tilde M(\zg))$ with triangulation $\tilde T(\zg)$.

Let $\mu=\mu_{j_s}\circ\mu_{j_{s-1}}\circ\ldots\circ\mu_{2}\circ\mu_1$
be a sequence of flips in $(S,M)$ such that the quadrilateral $V$
given by the arcs $\rho_1,\rho_2,\sigma_1,\sigma_2$ and $\zb'$ lies in
$\mu T$. Each $\mu_j$ flips an arc $\tau_j$ of a triangulation of
$(S,M)$. The fiber of $\tau_j$ in $(\tilde S,\tilde M)$ is of the form 
$\left((\tau_j)_g\right)_{g\in\Pi_1(S)}$. Let $\tilde \mu_j$ denote
the (infinite) sequence of flips in each $(\tau_j)_g$ in the fiber of
$\tau_j$ in $(\tilde S,\tilde M)$. Then
$\tilde\mu=\tilde\mu_{j_s}\circ\tilde\mu_{j_{s-1}}\circ\ldots\circ\tilde\mu_{2}\circ\tilde\mu_1$
is a sequence of flips in $(\tilde S,\tilde M)$ such that  every lift
$\tilde V$ of the quadrilateral $V$ lies in $\tilde\mu\tilde T$.

On the other hand, in $(\tilde S(\zg),\tilde M(\zg))$, we have the
lift $\overline V$ of $V$ given by the arcs  $\tilde \rho_1,\tilde
\rho_2,\tilde \sigma_1,\tilde \sigma_2$ and $\tilde \zb'$. There
exists a sequence $\overline\mu$ of flips in $(\tilde S(\zg),\tilde
M(\zg))$ such that $\overline V$ lies in $\overline \mu(\tilde
T(\zg))$.

Since  $(\tilde S(\zg),\tilde M(\zg))$  is a subsurface of $(\tilde
S,\tilde M)$, we can apply the same sequence of flips $\overline \mu$
to the triangulation $\tilde T$ of $(\tilde S,\tilde M)$ and see that
$\overline V$  lies in $\overline\mu\, \tilde T$. Clearly,
$\overline V$ is a lift of $V$ in $(\tilde S,\tilde M)$, and thus,
$\overline V$ also lies in $\tilde \mu\,\tilde T$.

$\pi:Q_{\tilde T}\to Q_T$ is a Galois covering of quivers, by Theorem
\ref{thm galois}.
Let $k$ denote the vertex in $\mu Q_T$ that corresponds to the arc
$\zb'$ and let $k_g$ denote the vertex in $\tilde \mu Q_{\tilde T}$ as
well as in $\overline \mu Q_{\tilde T}$ that correspond to the arc
$\tilde \zb'$ in $\overline V$.

For the purpose of this proof, we will use notations for the quiver $Q_{\tilde T(\zg)}$ that are induced from $Q_{\tilde T}$ via the inclusion $\tilde T(\zg)\subset \tilde T$. More precisely, the vertices of $Q_{\tilde T}$ are labeled as usual by tuples $(j,\zs)\in\{1,2,\ldots,2n\}\times\Pi_1(S)$ corresponding to arcs of $\tilde T$, if $1\le j\le n$, and to coefficients, if $n+1\le j\le 2n$;
and we use the same labels for the vertices of the subquiver $Q_{\tilde T(\zg)}$. Note that, by definition of $Q_{\tilde T(\zg)}$, there are precisely $d$ vertices that carry a label $(n+i, \zs)$ with $ 1\le i\le n$, and these correspond to those arcs $(i,\zs)$ in $\tilde T$ that are interior arcs in the triangulation $\tilde T(\zg) $ of the polygon $(\tilde S(\zg),\tilde M(\zg)) $.

Then, by Theorem \ref{thm galois}, the number of arrows $(n+i)\to k$
in $\mu Q_T$ is equal to the 
number of arrows $((n+i)_h)_{h\in \Pi_1(S)}\to k_g $ in $\tilde\mu
Q_{\tilde T}$. This last number is also equal to the
number of arrows  $((n+i)_h)_{h\in \Pi_1(S)}\to k_g$ in $\overline\mu
Q_{\tilde T}$, since
two sequences of mutations leading to the same quadrilateral will produce the same coefficients at the diagonals of the quadrilateral, because seeds are determined by their clusters, by \cite[Theorem 5.6]{FST}.
 In $Q_{\tilde T}$, the vertex $(n+i)_h$ is adjacent
only to the vertex $ i_{h}$, thus
if $(n+i)_h\to
k_g$ is an arrow in $\overline \mu Q_{\tilde T}$ then the sequence
$\overline \mu$ must flip the arc corresponding to the vertex 
$i_h$ at least once. By
\cite[pp. 40,41]{Mo}, we can conclude that  if $(n+i)_h\to
k_g$ is an arrow in $\overline \mu \,Q_{\tilde T}$, then the arc corresponding
to the vertex $i_h$ crosses $\zb'$ and hence, the arc corresponding to
$i_h$ is  an interior  arc of the polygon $(\tilde S(\zg),\tilde
M(\zg))$ and thus,  $(n+i)_h\to k_g$ is also an arrow  of
$\overline\mu Q_{\tilde T(\zg)}$.

Therefore, the number of arrows  $(n+i)\to
k$ in $\mu Q_{T}$ is equal to the number of arrows
$((n+i)_h)_{h\in\Pi_1(S)} \to
k_g$   in $\overline \mu Q_{\tilde T(\zg)}$, and consequently,
$y_\zb^+=\pi_*(\tilde y_{\tilde\zb}^+)$. The proof
of $y_\zb^-=\pi_*(\tilde y_{\tilde\zb}^-)$ is similar.
\qed
\end{pf}

\begin{lem}\label{lem pibar}
The covering map $\pi $ induces a bijection
$\pibar$ from the set of complete $(\tilde T(\zg),\tilde \zg)$-paths
in $\tilde S$ to  set of complete $( T, \zg)$-paths
in $ S$. 
 which sends a path
$\tilde \za=( \tilde \za_1,\tilde \za_2,\ldots,\tilde
 \za_{2d+1})$ to the path 
 $\pibar(\tilde \za)=
 ( \pi(\tilde \za_1),\pi(\tilde \za_2),\ldots,\pi(\tilde
 \za_{2d+1}))  $.
\end{lem}  
\begin{pf}
   
\emph{$\pibar $ is injective.} 
Suppose $\pibar(\tilde \za)=\pibar(\tilde \zb)$. Each $\tilde \za_i$
(respectively $\tilde \zb_i$) is the unique lift 
of $\pi(\tilde\za_i)$ (respectively $\pi(\tilde\zb_i)$) that starts at the
endpoint of $\tilde \za_{i-1}$ (respectively  $\tilde
\zb_{i-1}$), for $i=2,3,\ldots,d+1$. Since $\tilde \za_1$ and $\tilde\zb_1$ both start at
$\tilde s$ it follows that
$\tilde \za=\tilde\zb$ and $\pibar$ is injective.

\emph{$\pibar $ is surjective.} 
For every complete $(T,\zg)$-path $\za$ there is a  lift $\tilde
\za$ that starts at $\tilde s$. We have to show that $\tilde \za$ is a
complete $(\tilde T(\zg),(\tilde \zg))$-path. 
Since the crossing points of $\tilde \za$ and $\tilde \zg$ are $\tilde
p_1,\ldots,\tilde p_d$ in order,  $\tilde \za$ satisfies condition
(T1). Condition (T2) holds since $\tilde S(\zg) $ is simply connected. 
\qed
\end{pf}

\end{subsection}

\begin{subsection}{Cluster expansion in $\tilde {\mathcal
      A}$}\label{sect Atilde}
In this section, we will use Theorem \ref{thm main} for simply
connected surfaces to compute the cluster expansion of $\tilde
x_{\tilde \zg}$ in the cluster algebra $\tilde A$. Since $\tilde A$
does not have principal coefficients, we cannot use Theorem \ref{thm
  main} directly.
  
Let $x_{\tilde\zg}$ be the cluster expansion in the cluster algebra associated to $(\tilde S(\zg),\tilde M(\zg))$ with principal coefficient in the initial seed associated to the triangulation $\tilde T(\zg)$. By Theorem \ref{thm FZ4}, we have
\[ \tilde x_{\tilde \zg} = \frac {x_{\tilde \zg}\vert_{\tilde{\mathcal{F}}}(\tilde
  x_{i_1},\ldots,\tilde x_{i_d};\tilde y_{i_1},\ldots,\tilde y_{i_d})}
{x_{\tilde \zg}\vert_{\tilde{\mathbb{P}}}(1,1,\ldots,1,\tilde y_{i_1},\ldots,\tilde y_{i_d})}.
\]
Using Theorem \ref{thm main}, we get
\begin{equation}\label{eq 12}
\tilde x_{\tilde \zg} = 
    \frac {\sum_{\tilde \za}  x(\tilde \za)\vert_{\tilde{\mathcal{F}}}(\tilde
      x_{i_1},\ldots,\tilde x_{i_d}) 
         \  y(\tilde \za) \vert_{\tilde{\mathcal{F}}}(\tilde y_{i_1},\ldots,\tilde y_{i_d})}
{\sum_{\tilde \za}  y(\tilde \za)\vert_{\tilde{\mathbb{P}}}(\tilde y_{i_1},\ldots,\tilde y_{i_d})},
\end{equation} 
where both sums are over all complete $(\tilde T(\zg),\tilde
\zg)$-paths $\tilde \za$.
We compute the right hand  side of equation (\ref{eq 12}) in the following Lemma.

\begin{lem}\label{lem y}
Let $\tilde \za$ be a complete $(\tilde T(\zg),\tilde
\zg)$-path. Then

\begin{equation}\label{eq 13}
 x(\tilde \za))\vert_{\tilde{\mathcal{F}}}(\tilde
      x_{i_1},\ldots,\tilde x_{i_d})
=\frac {\prod_k \tilde x_{\tilde \za_k}}
{\tilde x_{i_1}\cdots\tilde x_{i_d}},
\end{equation}  
where the product is over all odd integers $k$ such that $\tilde
\za_k$ is an interior arc in $(\tilde S(\zg),\tilde M(\zg))$;

\begin{equation}\label{eq 14}
y(\tilde \za) \vert_{\tilde{\mathcal{F}}}(\tilde
          y_{i_1},\ldots,\tilde y_{i_d}) =
{\sum_{\tilde \zb}  y(\tilde \zb)\vert_{\tilde{\mathbb{P}}}(\tilde
          y_{i_1},\ldots,\tilde y_{i_d})} 
\left( \prod_\ell u_{i_\ell}
\right)
\left( \prod_k u_{\tilde \za_k}
\right)
\end{equation}  
where the sum is over all complete $(\tilde
T(\zg),\tilde \zg)$-paths $\tilde \zb$, the first product is over all integers
$\ell$ such that $\tilde \za_{2\ell}$ is $\tilde \zg$-oriented, and
the second product is over all odd integers $k$ such that $\tilde
\za_k$ is a boundary arc in  $(\tilde S(\zg),\tilde M(\zg))$.
\end{lem}

\begin{pf} Evaluating $ x(\tilde \za)$ in $\tilde
      x_{i_1},\ldots,\tilde x_{i_d}$ sets all the terms $\tilde
      x_{\tilde \za_k}$ with $\tilde \za_k$ a boundary arc in $(\tilde
      S(\zg),\tilde M(\zg))$ equal to $1$. Thus formula (\ref{eq 13})
      follows from the fact that all even arcs $\tilde \za_{2k}$ of
      $\tilde \za$ are interior arcs in  $(\tilde
      S(\zg),\tilde M(\zg))$.

      Formula (\ref{eq 14}) is a consequence of the following three
      claims.

  \emph{Claim 1: For $k=0,1,\ldots,d$, the exponent of $u_{[\zg_k]}$ in $ y(\tilde
  \za) \vert_{\tilde{\mathcal{F}}}(\tilde  y_{i_1},\ldots,\tilde y_{i_d})$ is
  $1$ if and only if $\tilde \za_{2k+1}=\tilde \tau_{[\zg_k]}$.}

      Proof of Claim 1. 
      By definition, 
      \[ y(\tilde \za) \vert_{\tilde{\mathcal{F}}}(\tilde  y_{i_1},\ldots,\tilde y_{i_d}) 
      = \prod_{k:\tilde \za_{2k} \textup { is $\tilde\zg$-oriented}} \tilde y_{i_k},
      \]
      and if $2\le k \le d-1$ then
      \[\tilde y_{i_k}=u_{i_k}\,u_{[\zg_k]}^{\zd_k}\,u_{[\zg_{k-1}]}^{\zd_{k-1}}, \textup{with}\]
      \[\zd_j=\left\{\begin{array}{ll}
      1&\textup{if $\tilde \tau_{[\zg_j]}$ follows $\tilde\tau_{i_k}$ in the clockwise orientation in $\tilde \zD_j$}    \\
     -1&\textup{if $\tilde \tau_{[\zg_j]}$ follows $\tilde\tau_{i_k}$ in the counter-clockwise orientation in $\tilde \zD_j$,}      
\end{array}\right.\]
for $j=k-1,k$.
Also,
 \[\tilde y_{i_1}=u_{i_1}\,u_{[\zg_0]}\,u_{[\zg_{-1}]}^{-1}, \textup{ and }
 \tilde y_{i_d}=u_{i_d}\,u_{[\zg_d]}\,u_{[\zg_{d+1}]}^{-1}. 
 \]
Therefore, the exponent of  $u_{[\zg_0]}$ is $1$ if and
      only if $\tilde \za_2$ is $\tilde \zg$-oriented, which is the
      case if and only if  $\tilde \za_1 =\tilde \tau_{[\zg_0]}$. This
      proves the case $k=0$. The exponent of  $u_{[\zg_d]}$ is $1$ if
      and
      only if $\tilde \za_{2d}$ is $\tilde \zg$-oriented, which is the
      case if and only if  $\tilde \za_{2d+1} =\tilde \tau_{[\zg_d]}$. This
      proves the case $k=d$.

      Now suppose that $0<k<d$.  The exponent of  $u_{[\zg_k]}$ is $1$
      if and only if, 
\begin{itemize} 
      \item[-] either $\tilde \za_{2k+2}(=\tilde \tau_{i_{k+1}})$ is $\tilde
      \zg$-oriented, $\tilde \za_{2k}(=\tilde \tau_{i_{k}})$ is not $\tilde \zg$-oriented,
      and the sequence \[\tilde \tau_{i_{k}},\tilde
      \tau{[\zg_k]},\tilde \tau_{i_k+1}\] corresponds to the
      counter-clockwise orientation of the triangle $\tilde \zD_k$,
      \item[-] or $\tilde \za_{2k}(=\tilde \tau_{i_{k}})$ is $\tilde
      \zg$-oriented, $\tilde \za_{2k+2}(=\tilde \tau_{i_{k+1}})$ is not $\tilde \zg$-oriented,
      and the sequence \[\tilde \tau_{i_{k}},\tilde
      \tau{[\zg_k]},\tilde \tau_{i_{k+1}}\] corresponds to the
      clockwise orientation of the triangle $\tilde \zD_k$.
\end{itemize}    
      In both cases, we have $\tilde \za_{2k+1}=\tilde\tau_{[\zg_k]}$,
      see Figure \ref{fig y}.

\begin{figure}
\input{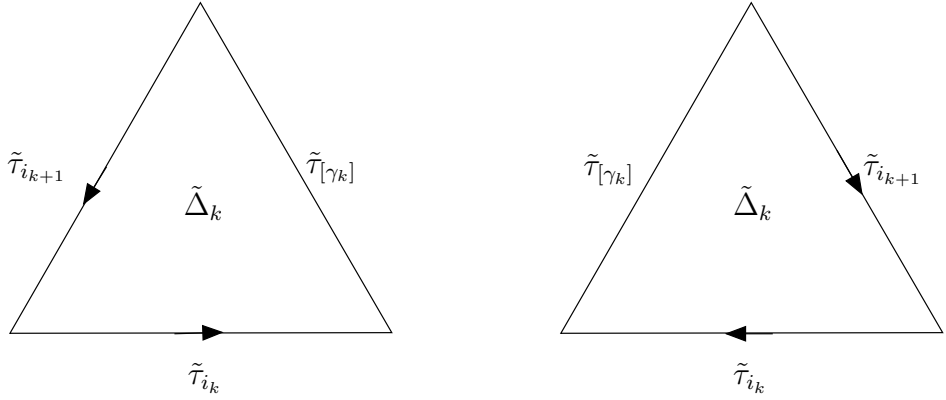}
\caption{The triangle $\tilde \zD_k$. On the left,  $\tilde \tau_{i_{k+1}}$
is $\tilde \zg$-oriented and $\tilde \tau_{i_{k}}$ is not; on the
right, $\tilde \tau_{i_k}$
is $\tilde \zg$-oriented and $\tilde \tau_{i_{k+1}}$ is not}\label{fig y}.

\end{figure}   
      
      Conversely, if $\tilde \za_{2k+1}=\tilde\tau_{[\zg_k]}$, then
      the segment 
\[(\tilde \za_{2k},\tilde \za_{2k+1},\tilde
 \za_{2k+2})=(\tilde \tau_{i_{k}},\tilde
 \tau_{[\zg_k]},\tilde \tau_{i_{k+1}})\] 
      goes around the boundary of $\tilde \zD_k$ exactly once. If this
      segment corresponds to the clockwise orientation of $\tilde
      \zD_k$ then $\tilde \za_{2k}$ is $\tilde \zg$-oriented and
      $\tilde \za_{2k+2}$ is not, and the exponent of $u_{[\zg_k]}$ is $1$; and if it corresponds to the
      counterclockwise orientation of $\tilde \zD_k$ then  $\tilde
      \za_{2k+2}$ is $\tilde \zg$-oriented and
      $\tilde \za_{2k}$ is not, and, again,the exponent of $u_{[\zg_k]}$ is $1$. This shows Claim 1.

   \emph{Claim 2: The exponent of $u_{[\zg_k]}$ in $ \tilde y(\tilde
  \za) \vert_{\tilde{\mathcal{F}}}(\tilde  y_{i_1},\ldots,\tilde y_{i_d})$ is
  $-1$ if and only if  
\begin{itemize}
\item[-] either $k=-1$ and $\tilde \za_{1}\ne \tilde\tau_{[\zg_{-1}]}$,
\item[-]  or $k=d+1$ and $\tilde \za_{2d+1}\ne \tilde \tau_{[\zg_{d+1}]}$.
\end{itemize}   }

      Proof of Claim 2. If $\tilde \za_1\ne \tilde\tau_{[\zg_{-1}]}$
      then  $\tilde \za_1 =\tilde\tau_{[\zg_0]}$, thus $\tilde \za_2 $
      is $\tilde\zg$-oriented and  the exponent of $u_{[\zg_{-1}]}$ is $-1$.
      Similarly, if  $\tilde \za_{2d+1}\ne \tilde\tau_{[\zg_{d+1}]}$
      then  $\tilde \za_{2d+1} =\tilde\tau_{[\zg_d]}$, thus $\tilde \za_{2d} $
      is $\tilde\zg$-oriented and  the exponent of $u_{[\zg_{d+1}]}$ is
      $-1$.

 Let us suppose now that $k\ne -1,d+1$ and
      $u_{[\zg_k]}$ has exponent equal to $-1$. Then 
\begin{itemize}
\item[-] either $\tilde
      \za_{2k}$ is $\tilde \zg$-oriented, $\tilde \za_{2k+2}$ is not
      and the sequence $\tilde \tau_{i_{k}},\tilde \tau_{[\zg_k]},
      \tilde \tau_{i_{k+1}}$ corresponds to the counter-clockwise orientation of
      the triangle $\tilde \zD_k$,
\item[-] or $\tilde
      \za_{2k+2}$ is $\tilde \zg$-oriented, $\tilde \za_{2k}$ is not
      and the sequence $\tilde \tau_{i_{k}},\tilde \tau_{[\zg_k]},
      \tilde \tau_{i_{k+1}}$ corresponds to the clockwise orientation of
      the triangle $\tilde \zD_k$.
\end{itemize}   
      Both cases are impossible. Indeed, in the first case, we would
      have \[(\tilde \za_{2k},\tilde \za_{2k+1},\tilde \za_{2k+2})
      =(\tilde \tau_{i_{k}},\tilde \tau_{i_{k+1}},\tilde
      \tau_{i_{k+1}}),\] contradicting the assumption that  $\tilde
      \za_{2k}$ is $\tilde \zg$-oriented and $\tilde \za_{2k+2}$ is
      not; and in the second case we would
      have \[(\tilde \za_{2k},\tilde \za_{2k+1},\tilde \za_{2k+2})
      =(\tilde \tau_{i_{k}},\tilde \tau_{i_{k+1}},\tilde
      \tau_{i_{k+1}}),\] contradicting the assumption that  $\tilde
      \za_{2k+2}$ is $\tilde \zg$-oriented and $\tilde \za_{2k}$ is
      not. This shows that we cannot have $k\ne -1,d+1$, and claim 2
      follows. 

  \emph{Claim 3: $\sum_{\tilde \zb}\tilde y(\tilde
  \zb)\vert_{\tilde{\mathbb{P}}}(\tilde
          y_{i_1},\ldots,\tilde y_{i_d}) = u_{[\zg_{-1}]}^{-1}
  u_{[\zg_{d+1}]}^{-1}$}.

      Proof of Claim 3. This follows from the definition of the
  addition $\oplus$ in $\tilde{\mathbb{P}}$, together with Claim 2 and  the fact
  that there always exist complete $(\tilde T(\zg),\tilde \zg)$-paths
  $\tilde \za$ and $\tilde \zb$ such that $\tilde
  \za_1=\tilde\tau_{[\zg_{-1}]}$ and $\tilde \zb_{2d+1}=\tilde \tau_{[\zg_{d+1}]}$.
\qed
\end{pf}  

Applying Lemma \ref{lem y} to equation (\ref{eq 12}), we get

\begin{prop}\label{prop y} 
Let $\tilde x_{\tilde \zg}$ be the cluster variable in $\tilde A$
corresponding an arc $\tilde \zg$ in $(\tilde S(\zg),\tilde M(\zg))$.
Then the cluster expansion of $\tilde x$ in the initial seed
associated to the triangulation $\tilde T(\zg)$ is given by
\[ \tilde x_{\tilde \zg} = \sum_{\tilde \za} \frac{\prod_j \tilde
  x_{\tilde \za_j} \prod_k u_{\tilde \za_k} \prod_\ell u_{i_\ell}}
  {\tilde x_{i_1},\ldots,\tilde x_{i_d}},
\]
where the sum is over all complete $(\tilde
T(\zg),\tilde \zg)$-paths $\tilde \za$, the first product is over all
odd integers $j$ such that $\tilde 
\za_j$ is an interior arc in $(\tilde S(\zg),\tilde M(\zg))$,
the second product is over all odd integers $k$ such that $\tilde
\za_k$ is a boundary arc in  $(\tilde S(\zg),\tilde M(\zg))$,
and the third product is over all integers
$\ell$ such that $\tilde \za_{2\ell}$ is $\tilde \zg$-oriented.
\end{prop}  
\end{subsection}

\begin{subsection}{Proof of Theorem \ref{thm main} in the general
    case}\label{sect general case} 
 In this subsection, we will prove
  Theorem \ref{thm main} for arbitrary unpunctured surfaces $S$.

Let $\zg$ be any arc in $(S,M)$, let $\tilde \zg$ be a lift of $\zg $
in $(\tilde S,\tilde M)$, and let $(\tilde S(\zg),\tilde
M(\zg))$ be the simply connected surface constructed in section
\ref{sect CA}. Let $\tilde A$ be the cluster algebra with combined
boundary and principal coefficients introduced in Definition \ref{def cl}. 
By Theorem \ref{thm cover}, we have $x_\zg=\pi_*(\tilde x_{\tilde \zg})$,
and from Proposition \ref{prop y} we get

\begin{equation}\label{eq 22}
   x_\zg = \sum_{\tilde \za} \frac{\prod_j \pi_*(\tilde
  x_{\tilde \za_j}) \prod_k \pi_*(u_{\tilde \za_k}) \prod_\ell
  \pi_*(u_{i_\ell}) }
  {\pi_*(\tilde x_{i_1},\ldots,\tilde x_{i_d})}, 
\end{equation}
where the sum is over all complete $(\tilde
T(\zg),\tilde \zg)$-paths $\tilde \za$, the first product is over all
odd integers $j$ such that $\tilde
\za_j$ is an interior arc in $(\tilde S(\zg),\tilde M(\zg))$,
the second product is over all odd integers $k$ such that $\tilde
\za_k$ is a boundary arc in  $(\tilde S(\zg),\tilde M(\zg))$,
and the third product is over all integers
$\ell$ such that $\tilde \za_{2\ell}$ is $\tilde \zg$-oriented.
By Theorem \ref{thm cover}, Lemma \ref{lem pibar} and Lemma \ref{lem or1}, it follows that 
\begin{equation}\label{eq 23}
   x_\zg = \sum_{\za} \frac{\prod_j x_{\za_j} }
  {x_{i_1},\ldots, x_{i_d}} \prod_\ell
  y_{i_\ell} ,
\end{equation}
where the sum is over all complete $(T,\zg)$-paths $\za$,  the first
  product is over all odd integers $j$ and the second product  is over
  all integers
$\ell$ such that $\za_{2\ell}$ is $ \zg$-oriented.
Thus
\[ x_\zg
=\sum_{\za} x(\za)\  y(\za),\]
where $\za$ runs over all complete $(T,\zg)$-paths in $(S,M)$, as required.
\qed

\end{subsection} 

\end{section}


\begin{section}{$F$-polynomials and $g$-vectors}\label{sect Fg}
In this section we study the $F$-polynomials and the $g$-vectors introduced in
\cite{FZ4}.

\begin{subsection}{$F$-polynomials}\label{sect F-polynomials}
By definition, the $F$-polynomial
$F_\zg$ is obtained from the Laurent polynomial $x_\zg$ (with
principal coefficients) by
substituting the cluster variables $x_1,x_2,\ldots,x_n$ by $1$.

\begin{thm}\label{thm F-polynomial}
Let $T$ be a triangulation of an unpunctured surface  $(S,M)$, and let
$\zg$ be an arc. Then the $F$-polynomial of $\zg$ is given by
\[ F_\zg = \sum_\za y(\za),\]
where the sum is over all complete $(T,\zg)$-paths in $(S,M)$.
\end{thm}  
\begin{pf} This follows immediately from Theorem \ref{thm main}.
\qed
\end{pf}

\begin{thm}\label{thm constant term}\cite[Conjecture 5.4]{FZ4}
Let $T$ be a triangulation of an unpunctured surface  $(S,M)$, and let
$\zg$ be an arc. Then the $F$-polynomial of $\zg$ has constant term $1$.
\end{thm}  
\begin{pf} We have to show that for any arc  $\zg$, there is
  precisely one complete $(T,\zg)$-path $\za^0$ such that none of its
  even arcs
  $\za_{2k} $ is $\zg$-oriented. As usual, let $d$ denote the
  number of crossings between $\zg$ and $T$. If $d=0$ then $\zg\in T$
  and there is nothing to prove. So suppose that $d\ge 1$ and let
  $p_1, p_2,\ldots,p_d$ be the crossing points of $\zg$ and $T$ in
  order of occurrence on $\zg$, and let $\tau_{i_j}\in T$ be such that
  $p_j$ lies on $\tau_{i_j}$.
Let $s$ be the starting point of $\zg$ and let $t$ be its
  endpoint. For $k=1,2,\ldots, d$, label the endpoints of $\tau_{i_k}$
  by $s_k,t_k$ in such a way that going along $\tau_{i_k}$ from $s_k$
  to $t_k$ is not $\zg$-oriented.
  Note that $s_k$ can be equal to $t_k$. 
  We have either $s_k=s_{k+1} $ or
  $t_k=t_{k+1}$. Also note that $i_k=i_{k'}$ does not necessarily
  imply that $s_k=s_{k'}$ since $\zg$ may cross $\tau_{i_k}$ several
  times and in both directions. Now define $\za$ as follows. 
 Let $\za_1=\tau_{[\zg_{-1}]}$ the unique arc in $T$ that
  is a side of the triangle $\zD_0$ and goes from $s$ to $s_{1}$. For
  $k=1,2,\ldots, d-1$, let $\za_{2k}$ be the arc $\tau_{i_k}$ running
  in the non $\zg$-oriented direction, and let $\za_{2k+1}$ be the arc
  $\tau_{i_k}$ running
  in the $\zg$-oriented direction if $s_{k+1}=s_k$, and let
  $\za_{2k+1}$ be the arc $\tau_{i_{k+1}}$ running in the
  $\zg$-oriented direction if $s_{k+1}\ne s_k$. 
Finally, let $\za_{2d}=\tau_{i_d}$ running in the non $\zg$-oriented
  direction, and let $\za_{2d+1}=\tau_{\zg_{[\zg_{d+1}]}}$ the unique arc in  in $T$ that is
  a side of the triangle $\zD_d$ and goes from $t_d$ to $t$.

Then $\za$ is a complete $(T,\zg)$-path and none of its even arcs is
$\zg$-oriented. 

To show uniqueness, it suffices to observe that specifying the even arcs
and their orientation, together with the homotopy condition (T2),
 completely determines the odd arcs and hence $\za$.
\qed
\end{pf}  

\begin{example} In Example \ref{ex A} the path $\za^0$ is the path 
\[(\tau_7,\tau_1,\tau_1,\tau_3,\tau_5,\tau_5,\tau_{11}),\] 
and in Example  \ref{ex Atilde} the path $\za^0$ is the path
\[(\tau_5,\tau_1,\tau_1,\tau_2,\tau_3,\tau_3,\tau_3,\tau_4,\tau_1,\tau_1,\tau_8 )\]
\end{example}
\begin{cor}\cite[Conjecture 5.5]{FZ4}
Each $F$-polynomial has a unique monomial of maximal
degree. Furthermore this monomial has coefficient $1$ and is divisible
by all the other occurring monomials.
\end{cor}  
\begin{pf} This follows from Theorem \ref{thm constant term} and
  \cite[Proposition 5.3]{FZ4}
\qed
\end{pf}  

\end{subsection}

\begin{subsection}{$g$-vectors}\label{sect g-vectors}
We keep the notation of section \ref{sect F-polynomials}.  Let $\tilde B_T$ be the $2n\times n$ matrix associated to the triangulation $T$ as in section \ref{sect quiver}, and let $B_T$ be the upper $n\times n$ part of $\tilde B_T$. It has been shown in \cite{FZ4} that, for any cluster variable $x_\zg$ in  $\mathcal{A}$, its Laurent expansion in the initial seed $(\mathbf{x}_T,\mathbf{y}_T,B_T)$ is homogeneous
 with respect to the grading given by
$\textup{deg}(x_i)=\mathbf{e}_i$ and
$\textup{deg}(y_i)=B_T\mathbf{e}_i$, where
$\mathbf{e}_i=(0,\ldots,0,1,0,\ldots,0) \in\mathbb{Z}^n$ with $1$ at
position $i$. By
definition, the \emph{$g$-vector} $g_\zg$ of a cluster variable $x_\zg$ is the 
degree of its Laurent expansion  with respect to this grading.

Because of Theorem \ref{thm constant term}, we also have $g_\zg=
\textup{deg}(x(\za^0))$, where $\za^0$ is the unique $(T,\zg)$-path
with $y(\za^0)=1$. As in the proof of Theorem \ref{thm
  constant term},  let $s_k$ be the startpoint of $\tau_{i_k}$ in
the non-$\zg$-oriented direction, $k=1,2,\ldots,d$.
Let 
\[
\begin{array}{rcl}
 I^-&=&\{i_k\mid k\in [1,d] \textup{ and } s_{k-1}=s_k\ne
s_{k+1}\} \\ 
I^+&=&\{i_k\mid k\in [2,d-1]\textup{ and } s_{k-1}\ne s_k=
s_{k+1}\}\cup\{[\zg_{-1}],[\zg_{d+1}] \},
\end{array}\]
where we use the convention that $s_0=s_1$ and $s_d \ne s_{d+1}$, so
that $I^-$ is well defined.
Then we have
\[ x(\za^0)= \prod_{h\in I^+} x_{h} \prod_{h\in I^-} x_{h}^{-1};\] 
recall that $x_{h}=1$ if $\tau_h$ is a boundary arc. 

Define $\mathbf{e}_{1},\mathbf{e}_{2} , \ldots, \mathbf{e}_{n}$ to be
the standard basis vectors of $\mathbb{Z}^n$, and let
$\mathbf{e}_{h}=(0,\ldots,0)$ if $\tau_h$ is a boundary arc.
We have proved the following theorem:
\begin{thm}\label{thm g-vectors}
Let $T$ be a triangulation of an unpunctured surface  $(S,M)$, and let
$\zg$ be an arc. Then the $g$-vector of $\zg$ is given by
\[g_\zg = \sum_{h\in I^+}\mathbf{e}_{h}-\sum_{h\in
  I^-}\mathbf{e}_{h} \]
\end{thm}  
\begin{example} In Example \ref{ex A}, we have 
\[I^-=\{3\} \quad I^+=\{7,12\} \quad g_\zg=(0,0,-1,0,0),\]
and in Example  \ref{ex Atilde} , we have
\[I^-=\{2,4\} \quad I^+=\{3,5,8\} \quad g_\zg=(0,-1,1,-1).\]
\end{example} 
\end{subsection}

Let us remark also the following useful fact.

\begin{lem}\label{lem degree}
Let $x_\zg=\sum_\za x(\za)\,y(\za)$ be the cluster expansion of
Theorem \ref{thm main} of a cluster variable $x_\zg$, and let
$\za,\zb$ be two complete $(T,\zg)$-paths.
Then 
\[x(\za)\,y(\za)=x(\zb)\,y(\zb) \ssi y(\za)=y(\zb) \]
\end{lem}  

\begin{pf} Suppose $y(\za)=y(\zb) $. Because of the homogeneity of the expansion proved in
  \cite{FZ4}, we have $\textup{deg}(x(\za))=\textup{deg}(x(\zb))$, and
  thus $x(\za)=x(\zb)$, since the degree of the  monomial $x(\za)$
  determines $x(\za)$ uniquely. Thus $x(\za)\,y(\za)=x(\zb)\,y(\zb) $.
The other implication is trivial.
\qed
\end{pf}  

\end{section}


\begin{section}{Arbitrary coefficients}\label{sect arbitrary}
Let $(S,M)$ be an unpunctured surface and $T$ a triangulation of
$(S,M)$. Denote by
$\mathcal{A}=\mathcal{A}(\mathbf{x},\mathbf{y},B_T)$ be the cluster
algebra with principal coefficients in the initial seed
$(\mathbf{x},\mathbf{y},B_T)$.

Let $(\hat{\mathbb{P}},\oplus,\cdot)$ be an arbitrary semifield,
let $\hat{\mathcal{A}}$ be a cluster algebra over the ground ring
$\mathbb{Z}\hat{\mathbb{P}}$ with an arbitrary coefficient system and initial seed
$(\mathbf{x},\hat{\mathbf{y}}, {B}_T)$, where
$\hat{\mathbf{y}}=\{\hat y_1,\ldots,\hat y_n\}$ is an initial
coefficient vector, $\hat {y}_i \in \hat{\mathbb{P}}$.

Let $\zg$ be an arc and $x_\zg$ and $\hat x_\zg$ the corresponding
cluster variables in $\mathcal{A}$ and $\hat{\mathcal{A}}$
respectively. Denote by $\hat{\mathcal{F}}$ the field of rational functions
$\mathbb{Q}\hat{\mathbb{P}}(x_1,\ldots,x_n)$. Then $\hat{\mathcal{F}}$ is the ambient field
of the cluster algebra $\hat{\mathcal{A}}$.

By Theorem \ref{thm FZ4}, we have
\[\hat x_\zg = \frac{x_\zg\vert_{\hat{\mathcal{F}}}(x_1,\ldots,x_n;\hat
  y_1,\ldots,\hat y_n)}
{F_\zg\vert_{\hat{\mathbb{P}}}(\hat y_1,\ldots,\hat y_n)},\]
where 
$x_\zg\vert_{\hat{\mathcal{F}}}(x_1,\ldots,x_n;\hat
  y_1,\ldots,\hat y_n)$ denotes the evaluation of the Laurent polynomial $x_\zg$ in
$  (x_1,\ldots,x_n;\hat   y_1,\ldots,\hat y_n)$ and
$F_\zg\vert_{\hat{\mathbb{P}}}(\hat y_1,\ldots,\hat y_n)$ denotes the
evaluation of the $F$-polynomial $F_\zg$ in $\hat
y_1,\ldots,\hat y_n$ in the semifield  $\hat{\mathbb{P}}$.

Using Theorem \ref{thm main}, we get

\begin{equation}\label{eq12}
\hat x_\zg =\frac{\sum_\za x(\za)\,\hat y(\za)}{\sum_\za \hat
  y(\za)\vert_{\hat{\mathbb{P}}}}, 
\end{equation}  
where $\hat y(\za)=y(\za)$ evaluated in $\hat y_1,\ldots,\hat y_n$; and
$ \hat  y(\za)\vert_{\hat{\mathbb{P}}}$ is $ \hat  y(\za)$ evaluated in
$\hat{\mathbb{P}}$, that is the addition $+$ is replaced by the addition
$\oplus$.

If $\hat{\mathcal{A}}$ is of geometric type,
i.e. $\hat{\mathbb{P}}=\textup{Trop}(u_1,\ldots,u_\ell)$ is a tropical semifield,
then $F_\zg\vert_{\hat{\mathbb{P}}}(\hat y_1,\ldots,\hat y_n)$ is equal to
$\prod_{i=1}^\ell u_i^{d(i)}$, where $d(i)$ is the 
minimum of all exponents of $u_i$
$F_\zg\vert(\hat y_1,\ldots,\hat y_n)$.
By Theorem \ref{thm constant term}, $F_\zg$ has constant term $1$, and thus $d(i)\le 0$.
In particular, $F_\zg\vert_{\hat{\mathbb{P}}}(\hat y_1,\ldots,\hat y_n)$ is
a monomial, and then Corollary \ref{cor pos} implies that the
right hand side of equation (\ref{eq12}) is a polynomial in
$x_1,\ldots,x_n$ with coefficients in $\mathbb{Z}_{\ge 0}\hat{\mathbb{P}}$. 
We have shown the following theorem known as the positivity conjecture.

\begin{thm}\label{thm pos}
Let $(S,M) $ be an
  unpunctured surface and let $\hat{\mathcal{A}}=\hat 
{\mathcal{A}}(\mathbf{x},\hat{\mathbf{y}},B)$ be a cluster algebra of
  geometric type associated to  some triangulation of $(S,M)$. Let $u$
  be any cluster variable   and let
\[u=\frac{f(x_1,\ldots,x_n)}{x_1^{d_1}\ldots x_n^{d_n}}\] be the
expansion of $u$ in the cluster $\mathbf{x}=\{x_1,\ldots,x_n\}$, where
$f\in\mathbb{Z}\hat{\mathbb{P}}[x_1,\ldots,x_n]$ .
Then  the coefficients of $f$ are non-negative integer linear
  combinations of elements in $\hat{\mathbb{P}}$.
\end{thm}  

\end{section}

\begin{section}{Further applications}
\begin{subsection}{Euler Poincar\'e Characteristics}\label{sect EPC}
In this section, we compare our cluster expansion formula to a formula
obtained recently by Fu and Keller \cite{FK} based on earlier formulas
in \cite{CC,CK,CK2,Palu}. Their formula is valid for cluster algebras
that admit a categorification by a 2-Calabi-Yau category, and, combining results of \cite{A} and \cite{ABCP,LF}, such a categorification exists in the case of cluster algebras associated to unpunctured surfaces. 

Particularly nice are the cases
where the surface is a disc or an annulus. In these cases 
the 2-Calabi-Yau category is the cluster category of type $A_{m-3}$, if
the surface is a disc with $m $ marked points, and of type  
$\tilde A_{p,q}$, if the surface is an annulus with $p$ marked points
on one boundary component and $q$ marked points on the other. 
Comparing the cluster expansions in these cases we shall obtain a
formula for the Euler-Poincar\'e characteristics of certain
Grassmannian varieties.

We recall some facts from cluster-tilting theory, see
\cite{BMRRT,BMR1,CCS1,ABS,CK2} for details.
Let $(S,M)$ be an unpunctured surface such that $S$ is either a disc
or an annulus.  As usual, denote by $n$ the number of interior arcs in any
triangulation of $(S,M)$.  If $S$ is an annulus, let  $p$ be the number of marked points
on one boundary component and $q$ the number of marked points on the
other. Define
\[ Q= \left\{
\begin{array}{ll}
\xymatrix@C=10pt@R=0pt{1\ar[r]&2\ar[r]&\cdots\ar[r]& n-1 \ar[r]& n
&
& \textup{if $S$ is a disc,}\\
\\
&2\ar[r]&\cdots\ar[r]&p\ar[rd]\\
1\ar[ru]\ar[rd]&&&&p+1&&
 \vspace{12pt} \textup{if $S$ is an annulus.}\\
&p+2\ar[r]&\cdots\ar[r]&p+q\ar[ru]\\
}
\end{array}  
\right.
\]
Let $A $ be the path algebra of $Q$ over an algebraically closed field
and let $\mathcal{C}$ be its cluster category. 

The indecomposable rigid objects of $\mathcal{C}$ are in bijection
with the cluster variables of the cluster algebra $\mathcal{A}$
associated to the surface $(S,M)$, hence, in bijection with the arcs
in $(S,M)$. Denote $X_\zg$ the indecomposable rigid object of
$\mathcal{C}$ that corresponds to the arc $\zg$. The clusters in $\mathcal{A}$, and
hence the triangulations of $(S,M)$, are in bijection with the
cluster-tilting objects in $\mathcal{C}$. Denote by
$T_{\mathcal{C}}=\{X_{\tau_1},X_{\tau_2},\ldots,X_{\tau_n}\}$ the cluster-
tilting object corresponding to the triangulation
$T=\{\tau_1,\tau_2,\ldots,\tau_n,\tau_{n+1},\ldots,\tau_{n+m}\}.$
The endomorphism algebra $\textup{End}_{\mathcal{C}}(T_{\mathcal{C}})$ of
this cluster-tilting object is called \emph{cluster-tilted
  algebra}. It has been shown in \cite{BMR1} that
$\Hom_{\mathcal{C}}(T_{\mathcal{C}},-)$ induces  an equivalence of categories  
$\mathcal{C}/\textup{add}\tau T_\mathcal{C}\to \textup{mod}\,
\textup{End}_{\mathcal{C}}(T_{\mathcal{C}})$. Let $M_\zg$ be the image of
$X_\zg$ under this equivalence.

Building on earlier results in \cite{CC,CK,CK2,Palu}, Fu and Keller
\cite{FK} proved a cluster expansion  formula for the cluster variable
$x_\zg$ which is of the form
\[ x_g= \sum_e \chi(Gr_e(M_\zg))\,x(e)\,y(e) ,\]
where the sum is over all dimension vectors $e=(e_1,e_2,\ldots,e_n)$,
$e_i\in\mathbb{Z}_{\ge 0}$; where $Gr_e(M_\zg)$ is the $e$-Grassmannian of
$M_\zg$, that is, the variety of submodules of dimension vector $e$;
where $\chi$ denotes the Euler-Poincar\'e characteristic; $x(e)$ is some
monomial rational function in the cluster variables $x_1,\ldots,x_n$,
and $y(e)=\prod_{i=1}^n y_i^{e_i}$. 

In particular, for any complete $(T,\zg)$-path $\za$, we have
$y(e)=y(\za)$ if and only if for each $k=1,2,\ldots,n$, the arc
$\tau_{i_k}$ appears exactly $e_{i_k}$ times as a $\zg$-oriented even
arc in $\za$. 
By applying Lemma \ref{lem degree}, we have proved the following theorem.

\begin{thm}\label{thm EPC} Let $M_\zg$ be a rigid indecomposable module of a
  cluster-tilted algebra of type $A$ or $\tilde A$. Then
$ \chi(Gr_e(M_\zg))$ is the number of complete $(T,\zg)$-paths $\za$
  such that  the arc
$\tau_{i_k}$ appears exactly $e_{i_k}$ times as a $\zg$-oriented even
arc in $\za$.
\end{thm}  

When $S$ is a disc this number is $0$ or $1$, but if $S$ is an annulus
it can be greater than $1$, see Example \ref{ex Atilde}. In
particular, we have the following Corollary.
\begin{cor}\label{cor EPC} $ \chi(Gr_e(M_\zg))$ is a non-negative integer.
\end{cor}  

\begin{rem} In the case where the cluster-tilted algebra is
  hereditary, Corollary \ref{cor EPC} was already proved in
  \cite{CR}. For general cluster-tilted algebras, it is new.
\end{rem} 

\end{subsection}

\begin{subsection}{Projective presentations}\label{sect projective
  presentations}
Let $T=\{\tau_1,\tau_2,\ldots,\tau_{n+m}\}$ be a triangulation of the
unpunctured surface $(S,M)$. Let $Q_T$ be the quiver defined in
section \ref{sect quiver}.
Let $A$ be the gentle algebra that the authors in \cite{ABCP}
associate to the triangulation $T$. The quiver $Q$ of the algebra $A$
is the full subquiver of $Q_T$ whose set of vertices is $\{1,2,\ldots,n\}$
and each  vertex  $j$ corresponds to an interior arc  $\tau_j$ of
$T$. For each
$j=1,2,\ldots,n,$ let $P(j)$ denote the indecomposable projective
$A$-module at the vertex $j$ of $Q$. 
Let $M_\zg$ be the
indecomposable $A$-module corresponding to the arc $\zg$. Its
dimension vector is given by the crossings of $\zg$ and the
triangulation $T$. 
Let $I^-$ and $I^+$ be the sets defined in section \ref{sect
  g-vectors}.

\begin{conj}\label{conj}
There is a minimal projective presentation of $A$-modules of the form
\[ \oplus_{h\in I^+} P(h) \to  \oplus_{h\in I^-} P(h) \to M_\zg \to 0.
\]
\end{conj}  

We prove the conjecture in the case where $S$ is a disc or an
annulus. With the notation of section \ref{sect EPC},
let 
\begin{equation}\label{eq 81}
 T^1 \to T^0 \stackrel{f}{\to} X_\zg \to T^1[1]
\end{equation}  
be a triangle in the cluster category $\mathcal{C}$, where $X_\zg$ is
the indecomposable rigid object corresponding to $\zg$, $T^1,T^0 \in
\textup{add}\, T_{\mathcal{C}}$ and $f$ is a right $\textup{add}\, T$-approximation.
Denote by $a_i^+ $ the multiplicity of the indecomposable $T_i$ in
$T^1$ and denote by $a_i^- $ the multiplicity of the indecomposable
$T_i$ in $T^0$. In other words, $T^1=\oplus_{i=1}^n a_i^+ T_i$ and
$T^0=\oplus_{i=1}^n a_i^- T_i$.
Applying $\Hom_{\mathcal{C}}(T_{\mathcal{C}},-)$ to the triangle (\ref{eq 81}), we
get a minimal projective presentation
\begin{equation}\label{eq 82}
 \oplus_{i=1}^n a_i^+ P(i) \to   \oplus_{i=1}^n a_i^- P(i) \to M_\zg
 \to 0
\end{equation}
in the module category of the cluster-tilted algebra. 
Following \cite{Palu}, we define the index $\textup{ind}\, M_\zg$ of $M_\zg$ to
be $(a_i^+-a_i^-)_i\in \mathbb{Z}^n.$ Then by \cite[Proposition 5.2]{FK},
we can relate $g$-vector and index by
\[g_\zg = - \textup{ind}\, M_\zg,\]
where the minus sign is due to our choice of orientation on
$(S,M)$. Thus, by Theorem \ref{thm g-vectors}, the projective
presentation (\ref{eq 82}) is of the form 
conjectured in \ref{conj}.

\end{subsection}

\end{section}


{} 

\vspace{1cm}
\noindent Ralf Schiffler\\
Department of Mathematics\\
University of Connecticut\\
Storrs, CT 06269-3009\\
Email: schiffler\@math.uconn.edu\\

\end{document}